\input epsf
\title{Remarks on AriasMarco-Sch\"uth's paper entitled:
``Local symmetry of harmonic spaces as determined by the spectra of
small geodesic spheres"
}

\newcommand{\cmt}[1]{\ifhmode\newline\fi{\sf *** \ \ #1 \\}}

\documentclass[11pt]{article}
\usepackage{latexsym}
\usepackage{epsfig}

\usepackage{amsmath}

\usepackage{amssymb}

\usepackage{exscale}

\usepackage{amsthm}
\usepackage{verbatim}

\newtheorem{theorem}{Theorem}[section]

\newtheorem{summary}[theorem]{Summary}

\def\:{\colon}

\long\def\onefigure#1#2{
\begin{figure*}[tbh]
\begin{center}
#1
\end{center}
\caption{#2}
\end{figure*}
} 

\newcommand{\iipefig}[1]  
{\smallskip\begin{center}\def\IPEfile{#1.ipe}\input{\IPEfile}\end{center}\smallskip}
\newcommand{\lipefig}[2]  
{\onefigure{\def\IPEfile{frh-#1.ipe}\input{\IPEfile}}{\label{f:#1} #2} }

\author
{Zolt\'an Imre Szab\'o
\thanks{Partially supported by NSF grant DMS-0604861}
\thanks{Lehman College of CUNY, Bronx, NY 10468,USA, and 
R\'enyi Institute, Budapest,
POBox 127, 1364 Hungary. 
}}
\date{}
\begin{document}

\maketitle
\begin{abstract}
The main goal in this paper is to point out that quantity 
$||\nabla R||^2(p)$ on a harmonic space can not be determined 
by the spectra of local geodesic spheres or balls, therefore the
main results of \cite{am-s} (quoted in the title) are wrong. 
My strong interest in this theorem is motivated by the
fact that it contradicts some of my isospectrality examples 
constructed on geodesic spheres and balls of certain harmonic manifolds.
The authors overlooked that 
the Lichnerowicz identity is not determined by the given spectral data,
and so is the final crucial equation obtained by eliminating 
with the Lichnerowicz identity. In short, 
the above theorem has falsely been established by  
spectrally undetermined identities which can not be computed (determined)
by the spectra of local geodesic spheres. More complicated 
spectrally undetermined functions cause the problems in 
case of local geodesic balls. I describe also a strong physical
argument which clearly explains why the manifolds appearing
in my examples are isospectral.  

It must be pointed out, however, that a very 
remarkable new idea, namely, the asymptotic expansion of the heat 
invariants $a_k(p,r)$ defined on geodesic spheres, $S_p(r)$, 
is introduced in the paper. 
It can be used for developing both geometric
uncertainty theory and global vs. local spectral investigations. 
Among my contributions to this developing field  
is the following statement: 
Average volumes, $\int vol(B_p(r))dp$ resp. $\int vol(S_p(r))dp$, 
of geodesic balls resp. spheres 
are generically not determined by the spectra of compact Riemann 
manifolds. This theorem interestingly contrasts the statement asserting 
that the volume of the whole manifold can be determined in terms of 
this global spectrum. In other words: Music written for local drums 
can not be played, in general, on the whole manifold.        
\end{abstract}

\section{A short history of harmonic manifolds.}  
Several equivalent definitions of harmonic manifolds are known in the 
literature. By one of them, they are the Riemannian manifolds yielding
the local mean value theorem of harmonic functions. By an other 
characterization, they are the manifolds where the density function
$\sqrt {\det (g_{ij})(q)}=\theta_p(q)/r^{n-1}$, defined on local normal 
coordinate systems
about an arbitrary point $p\in M^n$, is radial, that is, it depends just
on the distance $r(p,q)$. One can proof then that density function 
$\theta_p(r)$
does not depend on the reference point $p$ either.

The latter characterization clearly shows that all locally
two-point homogeneous manifolds are harmonic. 
The converse direction concerned,
Lichnerowicz (1946) showed that harmonic manifolds satisfying 
$dim M^n=n\leq 4$
are exactly the locally 2-point homogeneous manifolds. The question 
if this statement is true also in higher dimensions, became known as 
Lichnerowicz conjecture concerning harmonic manifolds. In 1990, this
problem was affirmatively answered on compact manifolds having finite
fundamental groups, by this author in \cite{sz1}. Due to the fact that all 
harmonic manifolds are Einstein, they have constant scalar curvature,
which is obviously strictly positive in the above compact cases. 
By Myers' and Cheeger-Gromoll's theorems, this affirmative answer 
extends also onto complete harmonic manifolds satisfying $Scal\geq 0$
\cite{sz2}.
However, extension onto non-compact manifolds satisfying $Scal<0$
defied every effort used by this author.

Then, in 1992, without knowing anything about the 
Lichnerowicz conjecture, Damek and Ricci \cite{dr} accidentally discovered 
infinitely many
non-compact and locally non-symmetric harmonic manifolds. Their statement
says that the natural invariant metrics on the solvable extensions, 
$SH^{(a,b)}_l$, of Heisenberg-type groups, $H^{(a,b)}_l$, are harmonic.
Such manifolds exist for any $l\geq 0$. Later we investigate families
defined by fixed values $(a+b)$ and $l$.
If $l\not =3\mod 4$, the metrics in such a family are isometric.
However, if $l\not =3\mod 4$, two metrics in the same family are locally
non-isometric, unless one of them is associated with $(a,b)$ and 
the other with $(b,a)$. Particularly interesting       
are the families defined for $l=3$. In fact, the metric on 
$SH^{(a+b,0)}_3\simeq SH^{(0,a+b)}_3$ is symmetric while the others 
are locally non-symmetric. 

Due to this discovery, the investigations were intensified 
both on compact and non-compact harmonic manifolds. Since the
Damek-Ricci examples do not have compact factors, the question arises
if the conjecture can be established on all compact manifolds and
not just on those satisfying $Scal\geq 0$.
By using entropy theory of geodesic flows, 
this question has been positively
answered on closed manifolds having negative sectional curvature, by G.
Besson, G. Courtois and S. Gallot \cite{bcg}, in 1995.
Since manifolds of strictly negative constant scalar 
curvature can have sectional curvature of changing 
sign, this question is still open on
closed manifolds belonging to this category.
In his very recent paper \cite{kn}, 
Knieper confirms the Lichnerowicz conjecture
for all compact harmonic manifolds without focal points or with
Gromov hyperbolic fundamental groups. These results strongly suggests
that the conjecture must be true on all manifolds that 
have compact Riemannian factors.
The Damek-Ricci spaces have also been very intensely investigated
in the past two decades. In this short list of achievements 
let it be mentioned just J. Heber's result \cite{h},
who, in 2006, proved that a homogeneous harmonic space is either
a model space or a Damek-Ricci space. The existence of inhomogeneous
harmonic manifolds is an  other important open question in this field.

By this author, spectral investigations 
of harmonic manifolds were initiated. Since then, this
field has been developed into several different 
directions.  This paper can focus mostly
on the investigations of this author. They
inevitably involve also the Heisenberg type 
groups whose solvable extensions are the Damek-Ricci spaces. 
Several different isospectrality examples have been 
constructed. The isospectral domains
arise always on the members of a family 
$H^{(a,b)}_l$ (resp. $SH^{(a,b)}_l$) defined by the same $(a+b)$ and $l$.  
The AM\&S-paper refers to the examples constructed in \cite{sz4,sz5},
where the isospectrality is confirmed for any pair of 
geodesic balls resp. spheres having the same 
radius $r$. 
The AM\&S-theorem contradicts the examples constructed on
the particular family $SH^{(a,b)}_3$, 
where the metric on
$SH^{(a+b,0)}_3\simeq SH^{(0,a+b)}_3$ is symmetric while the others in the
isospectrality family are locally non-symmetric. 
What makes this situation more serious is
that, in 2005, Hagen F\"urstenau \cite{f} discovered mathematical 
difficulties in the construction of the intertwining operator introduced
in \cite{sz4,sz5}. But, soon thereafter, this 
author did come up with the solution of this problem \cite{sz6,sz8}, 
which was publicly announced also at a conference held at CUNY, in
February, 2006 \cite{sz7}.   
Since then, I found also a strong physical argument explaining why
the isospectrality stated in \cite{sz4,sz5} must be true. Namely, it turns
out, that the Laplacian on the investigated manifolds can be identified
with the Hamilton operators of elementary particle systems and the 
isospectrality is equivalent to
the C-symmetry obeyed by these physical systems. 
This solution (together with explaining the difficulties in 
\cite{sz4,sz5})  
are described in the last section of this paper.  

\section{Technicalities on harmonic manifolds.}
Along a unit speed geodesic $c(r)$,
the Jacobian endomorphism field $A_{c(r)}(.)$ 
acting on the Jacobi subspace $V_{c(r)}\perp {\dot{c}(r)}=u(r)$ is defined 
by the equations
\begin{equation}
A^{\prime\prime}_{c(r)}+R_{\dot{c}(r)}A_{c(r)}=0,\quad A_{c(0)}=0, 
\quad A^\prime_{c(0)}=id,
\end{equation}
where $R_{\dot{c}(r)}(.)=R_u(.):=R(.,u)u$ is the Jacobian curvature
along $c(r)$. Its normalized version is denoted by 
$\mathbf A_{c(r)}=\frac 1{r}A_{c(r)}$.
The power series of  
invariants $Tr \mathbf A^k_{c(r)}$ can directly be computed by these
formulas. For instance, for $k=1$ we have:
\begin{eqnarray}
\mathbf A_{c_p(r)}=id +\frac 1{6}r^2R_{\dot{c}_p(0)}+
\frac 1{12}r^3R^\prime_{\dot{c}_p(0)}+ \frac {r^4}{5!}
(R_{\dot{c}_p(0)}^2+3R^{\prime\prime}_{\dot c_p(0)})
\label{boldA}
\\
+\frac {r^5}{6!}(4R^{\prime}_{\dot c_p(0)}
R_{\dot{c}_p(0)}+2R_{\dot{c}_p(0)}
R^{\prime}_{\dot c_p(0)}+4R^{\prime\prime\prime}_{\dot c_p(0)})+O(r^6),
\nonumber
\end{eqnarray}
where $c_p(r)$ indicates that it starts out from $p=c_p(0)$ and
$R^{\prime}_{\dot c_p(0)}(.)=R^{\prime}_{u_p}(.)
=(\nabla_{u_p}R)(.,u_p)u_p$,
etc.

Among these invariants 
the density $\theta_{c(0)}(c(r))=det\,A_{c(r)}$ plays the most important
role in this paper. Its
power expansion is usually
computed by Legendre's recursion formulas. Before providing the 
corresponding formulas computed by this direct 
method, we have to point out that some of the basic 
objects introduced here are different from those appearing 
in \cite{am-s} or \cite{be}. As opposed to these texts, the Laplacian is a
negative semi-definite operator and the
curvature is defined by 
$R(X,Y)Z=\nabla_X\nabla_YZ-\nabla_Y\nabla_XZ-\nabla_{[X,Y]}Z$
in all of my papers quoted later. In order to keep the exposition
simple, in this section we keep the AriasMarco-Sch\"uth-definition,
$R(x,y)z=-\nabla_x\nabla_yz+\nabla_y\nabla_xz+\nabla_{[x,y]}z$, 
for the curvature, but will return back to 
our definition in the rest part of this paper. There is
a difference also regarding $\theta$ in this paper which is associated 
with $det\,\mathbf A$ in \cite{am-s}. 
It should also be mentioned that some
of the constants appear in \cite{be} incorrectly not just in those
formulas which have been corrected in \cite{am-s} but also in some other
expressions which are quoted in this paper (these errors and their
effects on my computations were 
pointed out to me by D. Sch\"uth). In terms of these quantities we have:
 
\begin{theorem} (see \cite{be}, Chapter 6). If M is harmonic then 
there exist constants $C,H,L \in\mathbb R$
such that for all $p\in M$ and all $u\in T_pM$ with $|u| = 1$ we have:
\begin{eqnarray}
Tr(R_u) = C;\quad {\rm in \,\, particular}:
\label{C}
\\
Tr(R^{(k)}_u ) = 0,\quad {\rm for\,\, all\,\, 
derivatives\,\, of\,\, order} \,\, k \in \mathbb N.
\\
Tr(R_uR_u) = H;\quad {\rm in\,\, particular:}
\label{H}
\\
Tr(R_uR^\prime_u ) = 0\quad {\rm and}
\\
Tr(R_uR^{\prime\prime}_u ) =-Tr(R^\prime_u R^\prime_u ).
\\
Tr(32R_uR_uR_u - 9R^\prime_u R^\prime_u) = L.
\label{cond6}
\end{eqnarray}
\end{theorem}

The first five conditions imply that harmonic manifolds are Einstein, 
furthermore, the
curvature has constant norm $||R||^2=R_{abcd}R^{abcd}$.
Due to the Einstein property, these manifolds are analytic and the
unique radial density $\theta (r)$ is uniquely determined by the constants
appearing in the Taylor expansion. Such a density function 
does not determine the metric. See, for instance, that the above 
families of Damek-Ricci spaces share the very same density function 
$\theta (r)$. An other surprising feature of these examples is
that in families $SH^{(a,b)}_3$  
member $SH^{(a+b,0)}_3$ is symmetric while the others
are locally non-symmetric.

In this paper particular attention is paid to condition (\ref{cond6}).
It should be emphasized that this condition declares this 
curvature-expression to be 
constant just for this combination but it says nothing about 
the terms $R_uR_uR_u$ and
$R^\prime_u R^\prime_u$ which can very well be non-constant 
functions even along a
fixed geodesic.
The same density function can manifest itself very differently. For
instance, on $SH^{(a+b,0)}_3$, constant $L$ is incorporated completely
into $R_uR_uR_u$ but on the other members the same $L$ is distributed
both into $R_uR_uR_u$ and $R^\prime_u R^\prime_u$. Even
the proportion of the two terms depends on the members of a family.

By integrating these constants on the unit sphere formed by unit vectors
$u$ in the tangent space $T_p(M)$ one eliminates direction $u$ and
obtains the curvature conditions depending just on $p\in M$.
This computations result the following well known formulas:
\begin{eqnarray}
||R||^2=\frac{2n}{3}((n+2)H-C^2),
\label{curvnorm}
\\
32(nC^3+\frac 9{2}C||R||^2+\frac 7{2}\hat R-\mathring R)-27||\nabla R||^2=
n(n^2+6n+8)L,
\label{denseq}
\end{eqnarray}
where $\hat R$ and $\mathring R$ are defined in an orthonormal basis 
$\{e_1,\dots ,e_n\}$ of the
tangent space by the formulas:
\begin{equation}
\hat R =\sum_{i,j,p,q,r,s}R_{ijpq}R_{pqrs}R_{rsij},\quad
\mathring R =\sum_{i,j,p,q,r,s}R_{ipjq}R_{prqs}R_{risj}.
\end{equation}
 
Once again, the harmonicity conditions impose just this single equation
onto a tensor field having a large number of components. 
Thus, beyond this relation,
it provides no further information about the main players; 
$\hat R,\,\, \mathring R$, and $||\nabla R||^2$; in this formula. In other
words, these individual functions are out of touch of the harmonicity
conditions. Due to Lichnerowicz, there is an other independent equation:
\begin{equation}
\label{lich}
-\frac 1{2}\Delta ||R||^2=2C ||R||^2-\hat R-4\mathring R+||\nabla R||^2
\end{equation}   
which holds true not just on harmonic but on all Einstein manifolds.
On harmonic manifolds, relation $||R||^2=constant$ implies (see \cite{be},
Proposition 6.68):
\begin{eqnarray}
112\hat R-32\mathring R=27||\nabla R||^2+Q_1,\quad
\hat R+4\mathring R=||\nabla R||^2+Q_2,
\label{2eq}
\\
{\rm therefore},\quad\hat R=\frac {7}{24}||\nabla R||^2+Q_3, \quad 
\mathring R=\frac {17}{96}||\nabla R||^2+Q_4,
\label{elim}
\end{eqnarray}
where $Q_1,Q_2,Q_3,Q_4$ are constants which can be expressed in terms
of $C,H$, and $L$. (The second formula
on line (\ref{2eq}) appears with an incorrect factor 2 (before 
$||\nabla R||^2$) in \cite{be}.)    

It should be pointed out that there is a substantial difference 
between the equations (\ref{denseq}) and (\ref{lich}). The first
one is defined by integrating the constant function (\ref{cond6}) on the
unit sphere while no such natural integral representation exist for the 
Lichnerowicz identity. This identity is originally established by 
operations defined
just on $M$ which are nothing to do with 
curvature quantities such as $R_u^\prime R_u^\prime$;...; e. t. c., 
defined in terms of directions $u$. This detachment from functions 
defined by radial expansions strongly indicates that the Lichnerowicz
identity may cause serious difficulties in investigating spectrally
determined identities obtained by the radial expansion of the heat 
invariants $a_k(p,r)$. In fact, it will turn out that this identity is not
encoded into the spectra of local geodesic spheres, by any means, and
its application for reformulating spectrally determined identities
produces spectrally undetermined ones. This is why the identity
established in \cite{am-s} in the final step is not spectrally determined. 
Therefore, it can not prove the desired spectral determinacy stated
in the paper. There is also an other explanation for this problem
which does not use Lichnerowicz identity at all. 
Although it is defensible in case of geodesic spheres, it develops
into a major argument against the AM\&S-statement 
established by the spectra of local geodesic balls. 
The precise details are as follows.

\section{The AriasMarco-Sch\"uth local heat invariants.}

The first few heat invariants of a compact manifold $S$ are given by:
\begin{eqnarray}
a_0(S)=vol(S),\, a_1(S)= \frac1{6}\int_Sscal\, dvol_S,
\\
a_2(S)= \frac 1{360}\int_S(5scal^2-2||Ric||^2+2||R||^2)\, dvol_S.
\label{a2}
\end{eqnarray}
They are defined as coefficients in the asymptotic expansion
\begin{eqnarray}
\label{heatinv}
Tr(exp(t\Delta ))=\sum_{i=0}^\infty exp(\lambda_it)\sim
(4\pi t)^{-dim(S)/2}\sum_{k=0}^\infty a_k(S)t^k
\end{eqnarray}
where $t\downarrow 0$ and $\lambda_i\to -\infty$
denote the eigenvalues of Laplacian $\Delta$ with multiplicities.  

Arias-Marco and Sch\"uth furnished 
the spectral data of local geodesic
spheres by the asymptotic expansions of
radial functions $a_k(S_p(r))$, 
where $S_p(r)$ denotes the geodesic sphere of radius $r$
about $p$. The first few terms are explicitly computed for 
$a_2(S_p(r))$. The following theorem is a combination of 
Propositions 3.2 and 4.2 of \cite{am-s}, however, it is not a literal 
transcription of those statements. In order to prepare the arguments
discussed below, this transcription focuses also on the
eliminating steps (\ref{eq})-(\ref{reborn}). These are just briefly 
explained in \cite{am-s} without introducing any formulas.
 
\begin{theorem}\cite{am-s}
Consider an n-dimensional harmonic space $M$ with constants $C;H;L$, 
introduced above. Let $p\in M$ and let $u$ be a unit vector
in $T_p(M)$. Then
\begin{eqnarray}
||Ric^{S_p(r)}||^2_{exp(ru)}=\alpha_{-4}r^{-4}+\alpha_{-2}r^{-2}+\alpha_{0}
+\alpha_{2}r^{2}+O(r^3),
\\
{\rm and}\quad
||R^{S_p(r)}||^2_{exp(ru)}=\beta_{-4}r^{-4}+\beta_{-2}r^{-2}+\beta_{0}
+\beta_{2}r^{2}+O(r^3)
\end{eqnarray} 
for $r\downarrow 0$, where the coefficients $\alpha_{i}$ and $\beta_{i}$
for $i=-4,-2,2$ are constants depending only on $n,C$ and $H$. Moreover,
\begin{eqnarray}
\alpha_{2}(u)=\hat\alpha_{2}+\frac 1{16}Tr(R^\prime_u R^\prime_u)
\label{alpha2}
\\ 
{\rm and} \quad
\beta_{2}(u)=\hat\beta_{2}+\frac 4{9}
\sum_{i=1}^n Tr(R_u\circ R(e_i,.)R_ue_i),
\label{beta2}
\end{eqnarray}
where the coefficients $\hat\alpha_{2}$ and $\hat\beta_{2}$
are constants depending only on $n;C;H;L$ and $\{e_1,\dots ,e_n\}$
is an orthonormal basis of $T_p(M)$.

By integrating on the unit sphere $S_1(0_p)\subset T_p(M)$, whose volume
is denoted by $\omega_{n-1}$, we get:
\begin{eqnarray}
\int_{S_1(0_p)}Tr(R^\prime_u R^\prime_u)du=
\frac{3\omega_{n-1}}{n(n+2)(n+4)}||\nabla R||^2(p),
\label{Rprime}
\\
\int_{S_1(0_p)}\sum_{i=1}^n Tr(R_u\circ R(e_i,.)R_ue_i)du=
(nC^3+2\mathring R(p)-\frac 1{4}\hat R(p)) 
 \frac{\omega_{n-1}}{n(n+2)}.
\label{hatR}
\end{eqnarray}

By using (\ref{elim}) for eliminating terms $\mathring R(p)$ and 
$\hat R(p))$ from (\ref{hatR}) we get:
\begin{eqnarray}
2\mathring R(p)-\frac 1{4}\hat R(p)=\frac 9{32}||\nabla R||^2(p)
+Q(n,C,H,L),
\label{eq}
\end{eqnarray}
where $Q$ is a constant depending only on constants appearing in its 
argument, moreover:
\begin{eqnarray}
\overline\alpha_{2}:=\frac 1{\omega_{n-1}} \int_{S_1(0_p)}\alpha_{2}(u)du
=\tilde\alpha_{2}+ \frac{3}{16n(n+2)(n+4)}||\nabla R||^2(p),
\label{tildealpha}
\\
\overline\beta_{2}:=\frac 1{\omega_{n-1}} \int_{S_1(0_p)}\beta_{2}(u)du
=\tilde\beta_{2}+ \frac{9}{32n(n+2)}||\nabla R||^2(p),
\label{tildebeta}
\end{eqnarray}
where $\tilde\alpha_{2}$ and $\tilde\beta_{2}$ 
are constants depending only on $n,C,H$, and $L$.
Thus the scrutinized coefficient is reborn-ed in the new form: 
\begin{equation} 
\tilde\beta_{2}-\tilde\alpha_{2}+ \frac{9n-30}{32n(n+2)(n+4)}
||\nabla R||^2(p).
\label{reborn}
\end{equation}
\end{theorem} 

The computation proceeds by inserting this asymptotic expansion of
$a_2(S_p(r))$ into the integral formula 
$\int_{S_p(r)}a_2(S_p(r))\frac {dvol(S_p(r))}{vol(S_p(r))}$,
where the integration is taken with respect to the normalized density
$\frac {dvol(S_p(r))}{vol(S_p(r))}$. It is particularly important at
this step that the density is constant regarding $r$, that is, it appears
as the normalized density on the Euclidean unit sphere. Therefore,
in the asymptotic expansion of the whole integral formula this density
has no effect on the expansion. Also note that, due to the
uniquely determined density function $\theta (r)$ on a fixed harmonic
manifold (or for a whole family $SH^{(a,b)}_l$), 
the constants like $C$, $H$, $L$,...e. t. c. 
are spectrally determined values. After term by term integration
of the corresponding coefficient of the asymptotic 
expansion, the process is finished by eliminating $\mathring R(p)$
and $\hat R(p)$ by formulas (\ref{elim}). Then the final conclusion
is this: 
Since the coefficients of
the asymptotic expansion as well as the constant values are spectrally
determined, therefore also
$||\nabla R||^2(p)$ must be spectrally determined. 
This statement is considered also regarding the spectra 
of local geodesic balls. This case is discussed in the next section.

\section{Examining the AM\&S-proof.}
\subsection{Case of geodesic spheres.}
It must be emphasized again that the concerns brought up against 
this proof appear in the final eliminating step. All
preceding, highly non-trivial computations are performed with 
impressive carefulness. The key problem has already
been indicated at the introduction of the Lichnerowicz identity. 
According to those comments,
this equation does not have any natural integral representation
and its application cuts the computations completely off
from the spectral data. It seems to be that it is not encoded into 
the spectra of local geodesic spheres at all, thus, by adding it to
a local heat invariant, it may change it into a quantity which is not
local heat invariant any more.
Next, we show that exactly this is the case.
We will prove that this identity is really a spectrally 
undetermined object and when it is added, 
in the form of the above described 
linear combination with (\ref{denseq}), to the  
local heat invariant 
computed in the above theorem, 
then this elimination process changes this spectrally determined 
invariant to a spectrally undetermined curvature expression which
can not be produced as linear combination of ``true" local heat invariants 
obtained by direct radial expansions of several heat invariants 
$a_k(p,r)$. 
The conclusion with such a ``fake" local heat invariant
gives, of course, a false proof for the desired spectral determinacy
of $||\nabla R||$.    

In order to establish these statements, one should define 
the spectral data at a fixed point $p\in M$. It is defined by
the linear space of curvature expressions obtained by 
radial expansions of heat invariants $a_k(p,r)$. They can be computed
by integrating the coefficients defined by the radial expansions of
the heat invariant functions regarding $du_p$. The evaluations
of these integrals provide the curvature expressions which can be
written up as identities such that the right sides of these
identities are thought to be the integral values which are equal
to the curvature expressions obtained by
the AriasMarco-Sch\"uth expansions. On the left sides, the different 
curvature symbols can be
considered as independent objects (unknowns). These curvature terms
are defined on $M$, that is they depend on $p$ and not on the 
directions $u_p$. This identity interpretation
of these curvature expressions better indicates where these 
``true" local heat invariants come from. Thus the given
spectral data at $p$ is a linear space spanned by the identities 
obtained 
by the AriasMarco-Sch\"uth radial expansions of the complete set of 
heat invariants. 

In order to define
the complete linear space of local heat invariants at $p$,
one should introduce a naturally defined pre-Hilbert norm on this
linear space. It is defined by integrating
the products $F_1(u)F_2(u)$ of infinitesimal heat invariant functions
by means of $du$. By the topological closure associated
with this pre-Hilbert norm, 
one can define the {\it spectral-identity-space}
$SIS(p)$, which Hilbert space contains all spectrally determined 
identities attached to the given spectral data. 
The identities belonging to $SIS(p)$ are the only
ones which are determined by the given spectral data. 
Also note that this spectral data determines
just the identities obtained by integrating the infinitesimal heat 
invariant functions, but, it does not determines the pre-functions 
$F(u)$ from which these identities are derived. 

Because of this new interpretation, the
identities introduced so far must be rewritten in new appropriate forms.
This reformulation mostly concerns the constants appearing 
on the left sides of these
identities. Also these terms must appear as combinations
of curvature terms (unknowns) 
obtained by integrating the pre-functions. The values standing 
on the right sides are just formally symbolized by the integral 
$\int F(u)du_n$, where $du_n$ denotes normalized measure, 
whose actual evaluation provides the curvature terms, 
depending just on $p$, on the left sides of these equations. To this end,
consider the constants $C(u), H(u)$, and $L(u)$ defined in formulas
(\ref{C})-(\ref{cond6}) as functions of $u$ which are actually defined
by the left sides of those equations. Then $L$ on the right side of
(\ref{denseq}) must be substituted by $\int L(u)du_n$ while $C^3$ and
$C||R||^2$ on the left side should be replaced by the curvature 
expressions obtained by actual evaluations of integrals 
$\int C^3(u)du_n$ and $\frac{2n}{3}\int C(u)((n+2)H(u)-C^2(u))du_n$. These
integrals provide linear combinations of terms such as $Scal_M^3(p)$ and
$Scal_M(p)||R||^2(p)$. On non-Einstein manifolds these expressions are 
more complicated involving also the Ricci tensor and some other curvature
terms. Similar reformulation should be implemented in (\ref{lich}),
(\ref{alpha2}), (\ref{beta2}), (\ref{hatR})-(\ref{tildebeta}) and also
radial function $scal^2$ in (\ref{a2}) should be expanded. 

It is important to understand that this reformulation is not just a kind 
of fussiness. For instance, identities obtained by 
evaluating integrals like
$\int C^3(u)du_n$ and $\frac{2n}{3}\int C(u)((n+2)H(u)-C^2(u))du_n$
are not encoded into the spectra of local geodesic spheres, meaning
that they can not be obtained from the heat invariants $a_k(p,r)$ 
by radial expansions and linear combinations. 
That is, these curvature identities are not individually encoded into 
$SIS(p)$. The explanation for this phenomenon is that they can be
derived not from the original pre-functions defining the 
identities in $SIS(p)$ but from functions defined by
powering or multiplying the original heat invariant pre-functions.
They usually define curvature identities lying in the complement of
$SIS(p)$. In the following discussions such equations are
called {\it powered curvature identities}. 

Although identities associated with
$\int C^3(u)du_n$ or $\int C(u)H(u)du_n$ do not
show up in $SIS(p)$, yet, on harmonic manifolds, they are determined by
the spectra of local geodesic spheres. This statement is true for any
identity defined by $\int T_1(u)T_2(u)\dots T_k(u)du_n$,
where functions $T_i(u)$ are density-Taylor coefficients obtained by 
expanding $\theta$ regarding directions $u$.
Indeed, on harmonic manifolds,  
the Taylor coefficients of 
$\int\theta^k_p(u,r)du_n=(\int\theta_p(u,r)du_n)^k$ 
define spectrally determined identities, for any power $k\in\mathbb N$,
which are associated with integrals $\int T_1(u)T_2(u)\dots T_k(u)du_n$.  
In order to set up a spectral data which contains also these
powered density curvature identities, 
system $\{a_k(p,r)\}$ should be extended by the 
functions $\int\theta^k_p(u,r)du_n$ and extended spectral identity
space, $SIS^{ex}(p)$, should be defined by this complete 
extended function-system. This space incorporates all identities which
are determined by the spectra of local geodesic spheres of harmonic
manifolds.

It is interesting to see that how do these powered density functions
define the powered curvature terms. We demonstrate this by considering
the Taylor coefficient containing $L$ in the expansion of
$\theta$. It can be proved that, 
by increasing power $k$ by $1$, the
contribution to the previous function is 
a linear combination of terms $C^3,CH$, and $L$. Thus the
powered identities associated with these pre-functions can be
individually defined by powered density functions $\theta$, $\theta^2$,
and $\theta^3$. The AM\&S paper uses normalized density for the
expansions. In this paper we consider also the natural expansion
$a_k(p,r)$ without normalizing the density $\theta$ in the integral
formula. By the above discussion, this exchange of 
normalized density to the non-normalized one defines 
new coefficients for the terms
$L,C^3$, and $CH$ in the investigated Taylor coefficient 
associated with $L$. 
 
After these definitions we are able to explain why does the AM\&S-proof 
breaks down after eliminating with the Lichnerowicz identity. 
The ultimate problem is that the curvature expression standing on the left
side of this equation does not show up among the spectrally determined
curvature expressions. That is, this identity is not
encoded neither in $SIS(p)$ nor in $SIS^{ex}(p)$. This statement 
can be seen by proving that it is independent from both identity
spaces.

This observation is demonstrated, first, for the case satisfying
$||\nabla R||(p)\not =0$. Since it is independent from
the higher order equations expressed in terms of higher order 
curvature terms, it is enough to see that it is linearly 
independent from the lower-order spectrally determined identities 
introduced so-far. In case of $SIS(p)$,  
the problem reduces to
the independence of (\ref{lich}) from the system
consisting (\ref{denseq}) and the identity resulted by computations
(\ref{alpha2})-(\ref{hatR}). The desired independence can immediately 
be decided by considering only the main curvature expressions determined 
by terms $\hat R,\,\mathring R$, and $||\nabla R||^2$.
In case of $SIS^{ex}(p)$, the data contains also other identities,
but all new identities have main curvature expressions proportional
to the main term of (\ref{denseq}). Thus the independence is established
also in this case.

This independence
from the whole $SIS^{ex}(p)$ means that the curvature expression 
on the left side of the Lichnerowicz identity can not be produced by
linear combinations of curvature expressions obtained from
the expansions of $a_k(p,r)$ and $\theta^k(p,r)$.
It is obvious, that, the addition of a non-trivial 
constant-times of the Lichnerowicz identity to an identity belonging
to $SIS^{ex}(p)$ changes it to one which does not belong  $SIS^{ex}(p)$
any more. But exactly this happens during the considered elimination
which can be described such that a certain linear combination of the 
Lichnerowicz identity and (\ref{denseq}) is added to the 
spectrally determined identity resulted in (\ref{alpha2})-(\ref{hatR}).
Since (\ref{denseq}) is in $SIS^{ex}(p)$, it still keeps the identity
subjected to elimination in $SIS^{ex}(p)$. But the Lichnerowicz identity
results an identity which is certainly not in $SIS^{ex}(p)$ any more.
Therefore, this curvature identity (expressed in terms of 
$\int L(u)du_n, \int C^3(u)du_n,\int CH (u)du_n$, and $||\nabla R||^2$
but not containing neither $\hat R$ nor $\mathring R$) can not be
obtained by the radial expansions of functions $a_k(p,r)$ and
$\theta^k(p,r)$. Since all terms but $||\nabla R||^2$ can be obtained
by such expansions, thus $||\nabla R||^2$ is the only term which can
not be computed by such expansions. This proves that it is not 
a spectrally determined quantity.

In case of $||\nabla R||^2=0$, the above investigated identities
have only two main terms, $\hat R$ and $\mathring R$, thus the
independence of the Lichnerowicz identity from the other
two identity is not insured by the above arguments. If it is independent,
then one can imply in the same way that the elimination produces
a spectrally undetermined identity. But such an identity can not say
anything about the spectral determinacy of local symmetry. If it is
in $SIS^{ex}(p)$, then terms $\hat R$ and $\mathring R$ can be eliminated
from the investigated identity. After this operation only terms
associated with $\int C^3(u)du_n, \int C(u)H(u)du_n$ remain there. 
The only information provided by this step is that they are 
spectrally determined
quantities. But this is a trivial information which does not imply the
spectral determinacy of local symmetry, either.  

\subsection{Rudimentary vs. refined $SIS(p)$; Proper eliminations.} 
Before examining the 
proof of the AM\&S-statement in case of geodesic balls, 
we point out some of the rudimentary features of 
$SIS(p)$. As it turnes out, they can easily lead to misinterpretation of
the meaning of the spectral determinacy of the curvature
expressions defined by the given spectral data. 
Although the problems caused by these features can be fixed
in case of geodesic spheres, but they are not so
in case of geodesic balls. 

The primary problem about $SIS(p)$ is that
the curvature expressions obtained
by evaluating integrals $\int F_c(u)du_n$, where $F_c(u)$ is an 
$c^{th}$-order expansion coefficient defined by radial expansion of a
heat invariant $a_k(p,r)$, are considered without any
reference to the degree $c$. It is obvious that this data can not be
used for constructing functions depending on $r$. Such functions 
can be defined just by the complete expressions $\int F_c(u)du_nr^c$ 
involving $r^c$. This missing factor  
turnes the above elimination process to a very ambiguous operation. Indeed,
if a term from $\int F_c(u)du_n$ can be eliminated by a linear 
combination $\sum_iQ_i\int F_{c_i}(u)du_n$ where none of the
degrees $c_i$ is equal to $c$, then that term is still
present in $\int F_c(u)du_nr^c-\sum_iQ_i\int F_{c_i}(u)du_nr^{c_i}$. 
Furthermore, no linear
combination $A(p,r)=\sum_iP_ia_{k_i}(p,r)$ exists whose coefficient
of degree $c$ is equal to the curvature
expression resulting by the considered elimination.           

In order to explain this situation by a concrete example, suppose that
density coefficient $L$ appears in the coefficient 
$(\dots ,QL,\dots )_{n+1}r^{n+1}$ of $(n+1)$-degree 
in the natural expansion of a heat invariant $a_k(p,r)$ what we wish
to eliminate by the $L$ appearing in the coefficient of $(n+5)$-degree
in the natural expansion of the density heat invariant 
$\int\theta (p,r)du_n$. Then the coefficient $(\dots ,0,\dots )_{n+1}$ 
resulting by elimination can be decomposed as 
$(\dots ,0,\dots )_{n+1}=(\dots ,QL,\dots -QL)_{n+1}$,
meaning that it has $a_k(p,r)$- and 
$\int r^{-4}\theta (p,r)du_n$-components. Since the latter one is not a
heat invariant anymore, this elimination process involves new functions 
which
are not listed among the elements of the spectral data and produces
functions like $L(u)(r^{n+1}-r^{n+5})$ which still contain $L$
and which are not consistent with the requirement of radial 
constructibility of these functions. It is clear that all these problems
are caused by eliminating terms appearing in 
coefficients of distinct radial degrees  with each other. These arguments
can not be brought up for eliminations performed by coefficients
having equal degrees. Thus the only proper elimination in this situation
is the elimination of $QL_{n+1}$ appearing in 
$(\dots ,QL,\dots )_{n+1}$ by the  $L_{n+1}$ appearing in the 
corresponding coefficient defined by the natural expansion of
$\int r^{-4}\theta (p,r)du_n=r^{-4}\int\theta (p,r)$. Note that 
these coefficients are still
spectrally determined by the spectra of local geodesic spheres, because,
modulo certain multiplicative constants, just the degree defined by
$\int \theta (p,r)du_n$ must be lowered by $4$ in order to get the 
corresponding degree regarding $\int r^{-4}\theta (p,r)du_n$. Thus
the proper treatment of these problems requires just the extension
of the original spectral data $\{a_k(p,r)\}$ to the data 
$\{a_k(p,r),\int (r^{-4}\theta^k(p,r)du_n\}$.

A precise establishment of this refined data, 
$FineSIS(p)$, is as follows. In the first step the same type of expansion
should be chosen regarding each heat invariant $a_k(p,r)$, which can be
either the natural expansion or the normalized density expansion.
The degree, $c$, associated with a spectral identity 
$\int F_c(u)du_n=(Q\hat R+\dots )_c$ is indicated in the lower right
corners. They are the coefficients of $r^c$'s appearing in the 
complete radial expansion formulas.
For a fixed $c$, the $SIS_c(p)$ is generated by all heat invariant
coefficients of degree $c$. When this generator system is extended
by coefficients, of degree $c$, defined for all functions
$\int (r^{2d}\theta^k (p,r))du_n\}$, then the corresponding extended
coefficient space of degree $c$ is denoted by $SIS^{EXT}_c(p)$.
In order to avoid operations among coefficients of different degrees,
the total refined spectral identity spaces must be defined by
the direct sums $FineSIS(p)=\oplus_cSIS_c(p)$ resp. 
$FineSIS^{EXT}(p)=\oplus_cSIS^{EXT}_c(p)$. It is obvious that these 
spaces can be defined just for coefficients obtained from the
expansions of $a_k(p,r)$ resp. $\int (r^{2d}\theta^k (p,r))du_n\}$.

These refined spectral identity spaces put an end to the greatest
imperfection of the rudimentary $SIS(p)$'s which do not give any 
indication about the radial functions from which the elements of the
rudimentary $SIS(p)$ are derived. In case of $FineSIS(p)$ and 
$FineSIS^{EXT}(p)$,
however, any element $\sum_iQ_i\int F_c^{(i)}(u)du_n$ of
$SIS_c(p)$ or $SIS^{EXT}_c(p)$ determines the expansion term
$F_c(r)=\sum_iQ_i\int F_c^{(i)}(u)du_nr^c$ and there exist a function
$F(r)=\sum_cF_c(r)$ in the space spanned by functions 
$\{a_k(p,r)\}$ resp. $\{a_k(p,r),\int r^{2d}\theta^k (p,r))du_n\}$
whose $c^{th}$ expansion term is $F_c(r)$. Also the elimination
process appears in a well defined clear form. If the above linear 
combination $\sum_iQ_i\int F_c^{(i)}(u)du_n$ eliminates a term  
then $F(r)=\sum_cF_c(r)$ is such a globally defined function whose 
$c^{th}$ expansion coefficient is the curvature expression resulting by
the elimination. Eliminations which can be described by such globally 
defined functions are called {\it proper eliminations}. Only such
eliminations are acceptable to decide whether a term is spectrally
determined by the given data. (It is obvious that all those curvature
expressions are spectrally undetermined which can not be traced back
to the functions $\{a_k(p,r),\int r^{2d}\theta^k (p,r)du_n\}$ in this
way.) Eliminations performed by the elements of $SIS^{ex}(p)$
are called {\it rudimentary eliminations.}   
Let it be pointed out again that, by the above arguments, these
latter operations do not carry out the desired elimination, and, to the
curvature expression, which are the results of formal rudimentary 
elimination, no globally defined function exists which satisfies 
the above properties.

\subsection{Case of geodesic balls.}    

\medskip

\noindent{\bf A review of AM\&S-computations.}
In \cite{am-s}, the geodesic balls
are considered regarding both the 
Dirichlet and Neumann conditions. These computations start out from
the Branson-Gilkey formulas establishing the asymptotics of
the Laplacian on a manifold with boundary. After adopting these
generic formulas to harmonic manifolds, 
the Dirichlet coefficient $a_2^D(B)$ appears
in terms of the shape operator $\sigma_u(r)=(\nabla\nu )_{|T_{exp(ru)}M}$
(where $\nu$ is the normal vector field on $\partial B$) 
in the following form:
 \begin{eqnarray}
a_2^D(B)=\frac 1{360}\big(\int_BP_1^D(C^2,H)dvol_B +\int_{\partial B}
(P_2^D(..)+P_3(..))dvol_{\partial B}\big),  
\\
{\rm where}\quad
P_1^D(C^2,H)=5(nC)^2-2nC^2+\frac 4{3}n((n+2)H-C^2),
\\
P^D_2(..)=20nCTr(\sigma )-8CTr(\sigma )+16Tr(R_\nu\circ\sigma ),
\\
P^D_3(..)=\frac{40}{21}(Tr(\sigma ))^3-\frac{88}{7}Tr(\sigma )Tr(\sigma^2)
+\frac{320}{21}Tr(\sigma^3).
\end{eqnarray}

The radial expansions of these functions are established by means 
of the following formulas:
\begin{eqnarray}
Tr(\sigma )=(n-1)\frac 1{r}-\frac 1{3}Cr-\frac 1{45}Hr^3-
\frac 1{15120}Lr^5+O(r^7),
\\
Tr(\sigma^2)=(n-1)\frac 1{r^2}-\frac {2}{3}C-\frac 1{15}Hr^2+ 
\frac 1{3024}Lr^4+O(r^6),
\\
Tr(R_\nu\circ\sigma )=\dots +(-\frac 1{1440}L+
\frac 1{96}Tr(R_u^\prime R_u^\prime ))r^3+\dots ,
\\
Tr(\sigma^3)=\dots +(\frac 1{30240}L-\frac 1{96}
Tr(R_u^\prime R_u^\prime ))r^3+\dots . 
\end{eqnarray}
In the last two formulas only the $r^3$-coefficients of the corresponding 
series are indicated. Regarding the complete
functions $P_2^D$ and $P_3^D$, these coefficients
can be written in the form
\begin{equation}
[P_2^D]_3(CH,L,Tr(R_u^\prime R_u^\prime )),\,\, 
[P_3^D]_3(C^3,CH,L,Tr(R_u^\prime R_u^\prime )),
\end{equation} 
where also 
the arguments (on which these coefficients depend) are indicated.
 
Similar formulas are established also for $a_2^N(B)$. Although 
$P_1^D(C^2,H)=P^N(C^2,H)$ and $P^D_2=P^N_2$, the third one:
\begin{equation}
P^N_3(..)=\frac{40}{3}(Tr(\sigma ))^3+8Tr(\sigma )Tr(\sigma^2)+
\frac{32}{3}Tr(\sigma^3).
\end{equation}
is linearly independent from  $P^D_3$.
Thus, also 
$[P_3^N]_3(C^3,CH,L,Tr(R_u^\prime R_u^\prime ))$ is a new  
linear combination of the arguments. 

\medskip

\noindent{\bf Normalized-density expansions and powered terms on balls.}
For the radial expansion of this local heat invariant,
normalized density is used also in this part
of the paper. It is computed just regarding the term defined
on the boundary $\partial B$, in which case it means the expansion of
\begin{equation}
r\to\frac{1}{vol(S_r(p))}\int(P_2^D(..)+P_3(..))dvol_{S_r(p)}\big).
\end{equation}
In these discussions we need also the corresponding expansion of
the volume function. By the previous formula, this means the expansion of 
\begin{eqnarray}
r\to \frac{vol(B_r(p))}{vol(S_r(p))}
=\frac{1}{vol(S_r(p))}\int_0^r\int_{S_1(p)}\theta (u_p,\rho)du_pd\rho . 
\end{eqnarray} 
Note that this normalization concerns
just the density defined at the endpoint $\rho =r$ of radial integration
and it does not mean normalization regarding any $\rho$ between $0$ and
$r$. For explicit computations of the higher order powered terms appearing
in the formulas one must consider also the so called higher order 
k-volumes defined by 
$Vol_k(B_p(r))=\int_{B_p(r)}\Theta^k(p,\rho)dvol_{B_p(r)}$, where
$\Theta (r)=\det\mathbf A(c_p(r))=A_0+A_2r^2+A_4r^4+A_6r^6+\dots$.
The exact formulas, regarding both the natural and normalized-density
expansions of these functions, are described later in this section.
 
Before going into the details let it be mentioned in advance
that in case of geodesic spheres the spectrally
determined functions  $Vol_k(S_p(r))$ had to be extra added to the 
set $\{a_k(p,r)\}$ of
local heat invariants, because, they do not belong to them originally,
but, on harmonic manifolds, their spectral determinacy can be established 
by formulas such as $\int\theta^b(u_p,r)du_n=(\int\theta (u_p,r)du_n)^b$ 
and $\int\Theta^b(u_p,r)du_n=(\int\Theta (u_p,r)du_n)^b$.
This spectral determinacy can not be established for functions 
$Vol_k(B_p(r))$, however. Indeed,  
for $k>0$, these functions are nothing to do with the heat invariants
of local geodesic balls, furthermore, there does not exist a 
function-relation
$Vol_k(B_p(r))=F_k(vol(B_r(p))$ expressing the k-volume $Vol_k(B_p(r))$
in terms of the standard volume $vol(B_p(r))=Vol_0(B_p(r))$, for 
all $r$. Thus
neither the k-volumes nor the proper powered terms (obtained by using
higher order k-volume-functions for eliminations) are determined by the
spectra of local geodesic balls. Therefore this spectral data can not
be extended by these spectrally undetermined functions.   

\medskip

\noindent{\bf Rudimentary vs. proper eliminations on balls.}
After this detour we turn back to the AM\&S-proof,
where the desired spectral determinacy of $||\nabla R||^2$ is  
concluded just by a briefly explained argument. According to this 
explanation, just term $L$  
appearing in the above $r^3$-coeffi\-ci\-ents has to be eliminated by 
the $L$ appearing in the corresponding 
coefficient of the radial expansion of the first heat invariant 
$a_0^D(p,r)=vol(B_r(p))=a_0^N(p,r)$. 
(The elimination of this $L$ is really
necessary because it involves $Tr(\nabla_uR_u\nabla_uR_u)$.) 
However, a deeper insight into the
details shows that this briefly 
explained step covers up numerous problems
which discredit the AM\&S-proof also in this case. One of them is
indicated above. Even if this elimination were appropriate, there are
still higher order spectrally undetermined powered terms left behind
because of which term  $||\nabla R||^2$ can not be spectrally determined.
But in reality there are much more deeper reasons why this quantity is
spectrally undetermined.

First of all, the computations below show that the $L$ appears
with different degrees regarding the two functions  $a_0(p,r)$ and
$a_2(p,r)$, with respect to both type of expansions. The degree computed
by $a_2(p,r)$ is always lower by $2$ than the degree computed 
by $a_0(p,r)$.
This means that the proposed elimination is a typical rudimentary
operation which leads to several contradictory
situations. Remember that, in order to keep connection with the radial
functions from which these coefficients are defined, this operation 
must be well defined in terms of these functions. Since radial functions
$Q_1Lr^c$ and  $Q_2Lr^{c+4}$ can not cancel out each other, the desired
elimination of $L$ does not take place, furthermore, the coefficient 
defined by the formal rudimentary elimination does not appear as a  
coefficient of a spectrally determined function, which,
by definition, must be a linear combination of functions  
$\{a_k(p,r)\}$. In the top of these problems, also
spectrally undetermined higher order powered terms are present 
in both functions. Since they
can not cancel out each other, the 
proposed rudimentary operation is totally inappropriate for the
desired elimination.

Proper eliminations can be implemented by functions 
\begin{equation}
\int_{B_p(r)}\rho^{-4}\Theta^k(p,\rho)dvol_{B_p(r)},\quad {\rm where}
\quad k\geq 0,
\end{equation} 
but this adjustment changes even $vol({B_p(r)})$ into a spectrally
undetermined function and it further worsens the spectral indeterminacy
of the k-volumes of balls. Remember that, in case of geodesic spheres,
the spectral determinacy of 
$\int r^s\theta^k(u_p,r)du_n=r^s\int(\theta(u_p,r))^kdu_n$
is due to the indicated equation. Since $\rho^{-4}$ is inside of an
integration, such reconnection to the volume does not exist in case
of geodesic balls. (By using multiple partial integrations, term
$\rho^{n-5}$ can be retraced to $\rho^{n-1}$. The complicated other
terms produced by this process show then the impossibility of the
reconnection of these adjusted functions to the volume.) In other
words, all these functions are spectrally undetermined, together with
all those terms which can be eliminated by these functions in a proper
elimination process. This means that the $L$, which can not 
properly be eliminated from $a_2(p,r)$ by a spectrally determined 
function, but just by a spectrally undetermined one, is spectrally 
determined in  $a_0(p,r)$ whereas it is spectrally undetermined in  
$a_2(p,r)$. That is, the spectral determinacy of a term such as $L$
is not an absolute but a relative concept. This property depends 
on the degree of the expansion term in which it shows up. This
paradoxical statement points to the deep sitting ambiguity present in
rudimentary spectral data defined for geodesic balls. It clearly
explains also the error made in this part of \cite{am-s}. By the
proposed elimination they tried to cancel out a spectrally undetermined
quantity by a spectrally determined one. As it is explained above,
this is a totally inappropriate operation from several point of view.     

\medskip

\noindent{\bf Explicit expansion formulas.}
In this section those expansions are explicitly computed
which are mentioned in the above discussions. 
Expansion of $vol(B_r(p))/vol(S_r(p))$ appears in the form: 
\begin{equation}
\frac{vol(B_r(p))}{vol(S_r(p))}= \frac{vol(B_r(p))}
{vol^\prime (B_r(p))}=r(V_0+V_2r^2+V_4r^4+V_6r^6+\dots ).
\end{equation}
The details below show that this expansion does not provide
$L$ in a pure form. Instead, it is accompanied with higher order 
powered terms.
In fact, a simple calculation shows that coefficients $V_k$ can be
determined by the coefficients of expansions 
\begin{eqnarray}
\Theta (r)=\det\mathbf A(c_p(r))=A_0+A_2r^2+A_4r^4+A_6r^6+\dots ,
\\
vol(B_r(p))=r^n(\frac 1{n}A_0+\frac 1{n+2}A_2r^2+\frac 1{n+4}A_4r^4+
\frac 1{n+6}A_6r^6+\dots )
\label{vol(B)}
\end{eqnarray}
by the recursive formula
\begin{equation}
V_k=\frac 1{A_0}(\frac 1{n+k}A_k-A_2V_{k-2}-A_4V_{k-4}-\dots -
A_kV_{0}),
\end{equation}
where $k$ denotes even numbers.
Thus (\ref{denseq}) really does not arise in a pure form in this 
expansion. Notice that, besides $L$, powered terms like $C^3$ and $CH$ 
are also present in these local heat invariants. Thus the proposed
rudimentary elimination contributes new proper powered terms to the
ones already present there. 
Since the degrees defined by the considered two functions are different,
these complications create new problems for the proposed rudimentary
elimination. This phenomenon more emphatically underlines the
inappropriateness of the proposed rudimentary elimination. 

Also note that the expansion by normalized density represents these 
volumes $Vol_k(B_r(p))$, where $k\geq 0$, by series of the form 
$r(V_{k0}+V_{k2}r^2+V_{k4}r^4+V_{k6}r^6+\dots )$ and such $O(r)$-series 
appears just regarding the first function $\int_BP_1^D(C^2,H)dvol_B$
of $a_2(p,r)$.
The other two terms define $O(\frac 1{r})$- and 
$O(\frac 1{r^3})$-series. The $r^3$-component regarding the
first function with respect to normalized-density expansion is a constant
time of $CP_1(C^2,H)$, while this function is  $LP_1(C^2,H)$ regarding
the $r^7$-component.
 
The above arguments can easily be established also for 
natural expansions. Then
function $vol(B_p(r))$ appears in the form described in 
(\ref{vol(B)}). That is, no proper powered terms appear in this case
which makes these computations more 
transparent than they were in case of normalized density expansions. 
The natural expansion of the functions defined on the boundary
$\partial B$
can be established by expanding also $\theta$ and multiplying
this series with the AM\&S-series. Since functions $P_2^D$ and $P_3^D$
have singularities of order $O(\frac 1{r})$ and $O(\frac 1{r^2})$
respectively, furthermore, $\int_BP_2^Ddvol_B$ involves 
integration regarding $dr$, this Taylor series can be written in
the form
$
\sum_k(Q^D_{1(k-4)}+Q^D_{2(k-2)}+Q^D_{3k})r^{n-4+k}
$,
where $Q^D_{ik}$ is associated with the product of 
$P^D_{i}$ and $\theta$. That is: $Q^D_{1}=r^n(Q^D_{10}+Q^D_{12}r^2+..$,
$Q^D_{2}=r^{n-2}(Q^D_{20}+Q^D_{22}r^2+..$, and 
$Q^D_{3}=r^{n-4}(Q^D_{30}+Q^D_{32}r^2+..$. 
The above described $r^3$-coefficients appear then
in terms of $Q^D_{12}(C^3,CH)$, 
$Q^D_{24}(CH,L,Tr(R_u^\prime R_u^\prime ))$,
and $Q^D_{36}(C^3,CH,L,Tr(R_u^\prime R_u^\prime ))$ where,
comparing with those in $[P^D_i]_3$, the
arguments show up in completely new linear combinations. 
Indeed, the product of functions $P^D_i$ with
$\Theta$ provides different contributions to the arguments of 
the individual terms labeled by $i=1,2,3$.
Unlike in case of normalized density expansions, 
no powered terms show up in the natural expansion formulas. But,
because of the different degrees, the proposed rudimentary elimination 
still remains a completely inappropriate operation also regarding this
expansion.

\medskip

\noindent{\bf Other explanations for the spectral indeterminacy 
of $||\nabla R||^2$.} In this section we describe two new 
interrelated constructions which also prove the investigated
spectral indeterminacy.

{\bf (1)} In the previous constructions the required elimination
of $L$ and the other higher order density terms can 
properly be performed just by adjusted
volume functions. More precisely, this operation can be completed by the 
adjusted standard volume, $Vol^{(adj,\rho^{-4})}_0$ and
two higher order k-volumes, $Vol^{(adj,\rho^{-4})}_1$ and 
$Vol^{(adj,\rho^{-4})}_2$. 
One can go in the other way around, by adjusting
$a_2(p,r)$ to the volume functions $Vol_0,\, Vol_1$, and $Vol_2$.
This adjusted function is then
 \begin{eqnarray}
a_2^{(adj,1,r^2,r^4)}(B)=\frac 1{360}\big(\int_BP_1^D(C^2,H)dvol_B +
\\
+\int_{\partial B}
(r^2P_2^D(..)+r^4P_3(..))dvol_{\partial B}\big),  
\end{eqnarray}
where functions $P_1,\, P_2$, and $P_3$ are multiplied with 
$1,\, r^2$, and $r^4$ respectively. This adjustment transforms the local
heat invariant function to a spectrally undetermined function, which, 
like the volume function
$vol(B(r))$, is of order $O(r)$ in case of normalized density 
expansions and of order $O(r^n)$ in case of natural expansions. 
In each cases the $L$ and the
third order power terms defined on the boundary $\partial B$ by the
last two integrals involving $P_2$ and $P_3$ appear with the same degree
as $L$ does regarding  $vol(B(r))$. Also notice that the first integral
$\int_BP_1^D(C^2,H)dvol_B$ defines term $LP_1(C^2,H)$ regarding this 
degree. Term $Tr(\nabla_uR_u\nabla_uR_u)$ present in this coefficient
of $a_2^{(adj,1,r^2,r^4)}(B)$ can be obtained then by eliminating 
all these density terms by means of k-volume functions 
$Vol_0,Vol_1,Vol_2,Vol_3,\dots$. Because of term $LP_1(C^2,H)$,
this process can be completed by using more higher order volume functions
than those used in the previous cases. Among these functions, 
which set includes also
$a_2^{(adj,1,r^2,r^4)}(B)$, only $Vol_0$ is spectrally determined. Thus
also this representation defines 
$Tr(\nabla_uR_u\nabla_uR_u)$ as a spectrally
undetermined function. The same statement can be established by the 
adjusted
functions $a_2^{(adj,q_1,q_2r^2,q_3r^4)}(B)$, where $q_1,q_2,q_3$ are
arbitrary constants.

{\bf (2)} The above construction strongly suggests that the
spectral determinacy of $||\nabla R||^2$ can be established by splitting
$a_2(B_p(r))$ into 3 different functions. 
The decomposition of $a_2(p,r)$ into three different parts is one
of the most important characteristic features of these functions.  
The components in the coefficient
$(Q_{12}+Q_{24}+Q_{36})$ form
a linearly independent system of linear functions depending on the
variables $C^3,CH,L$, and $Tr(R_u^\prime R_u^\prime )$. System 
$\{L, Q_{12}, Q_{24}, Q_{36}\}$ consists of for independent linear
functions of these four variables, which appear as the 
corresponding Taylor
coefficients of functions  
$\{a_0(p,r), Q_{1}(p,r), Q_{2}(p,r), Q_{3}(p,r)\}$, respectively.
By the given spectral data only 2 functions: 
\begin{equation}
\{a_0(p,r), \sum_k (Q_{1(k-4)}(p) +Q_{2(k-2)}(p)+ 
Q_{3k}(p))r^{n-4+k}\}
\end{equation}
are determined by which no splitting of the second function can be
established. It is the minimal data which provides 
the desired spectral determinacy of $||\nabla R||^2$.

The same
statements are true also for the spectral identities obtained by 
evaluating the integrals.  
That is, identity
\begin{eqnarray}
\int (Q_{12}+Q_{24}+Q_{36})du=(\tilde Q_{12}(Scal^3, Scal||R||^2)+
\\
+\tilde Q_{24}(Scal^3, Scal||R||^2,\frac 7{2}\hat R-\mathring R,
||\nabla R||^2)+
\nonumber
\\
+\tilde Q_{36}(Scal^3, Scal||R||^2,
\frac 7{2}\hat R-\mathring R,||\nabla R||^2))
\nonumber
\end{eqnarray}
is the sum of 3 linearly independent identities which define a 
spectrally determined identity just in this linear combination.
These component identities together with spectral
identity $\tilde L(Scal^3, Scal||R||^2,\frac 7{2}\hat R-\mathring R,
||\nabla R||^2)$ also determine a linearly independent system. The
latter one is available now by the expansion of $a_0^D$.
The desired spectral determinacy of 
can be pointed out just by a drastically changed $SIS(p)$ which, 
along with 
$\tilde L(Scal^3, Scal||R||^2,\frac 7{2}\hat R-\mathring R,
||\nabla R||^2)$, separately contains all above 3 spectral identities.
That is, in case of geodesic balls, the split of the heat invariants
into three parts provides a new spectral data by which the desired 
spectral determinacy can be established. Note that this is the minimal
data by which the investigated problem can be positively answered.
All attempts must fail unless they are able to establish this division
of the heat invariant into 3 separate invariants. But this split is
certainly not encoded into the spectra of local geodesic balls, thus,
all these attempts must fail altogether. 

\medskip

\noindent{\bf Unified spectral data.} 
Note that, regarding the two boundary conditions, 
the first two components in the above decompositions are the same, 
thus the difference eliminates both ones, leaving
behind just $P^D_3-P^N_3$.  
Therefore, it seems to be possible to squeeze out some kind of positive
answer regarding the spectral determinacy of $||\nabla R||^2$ by
eliminating $\tilde L$, $Scal^3$, and $Scal||R||^2$ by using 
the two conditions together such that the unified data 
$\{a_k^D(p,r),a_k^N(p,r)\}$
is extended  by functions $Vol_k(p,r)$. 
But also this elimination is just rudimentary which does not 
work out by the very same reason explained earlier. In all of these 
cases proper elimination can be defined just by the refined spectral
data $FineSIS_B(p)=\oplus_kSIS_{B,k}(p)$ which is the direct sum of 
$SIS_{B,k}(p)$'s defined for radial degrees $k=0,1,2,3,\dots$. 
This definition excludes eliminations of terms with other ones
which have different radial degrees. The positive conclusion can
not be squeezed out even if the data is enlarged by the spectra of 
local geodesic spheres. In this case functions $\int\theta^k(p,r)du_p$
exhibit contradictory properties regarding these two different type of
spectra. Actually they are incompatible which is indicated
by the distinct factors $(4\pi t)^{-n/2}$ and $(4\pi t)^{-dim(S)/2}$ 
appearing in the asymptotic expansions of
the partition functions, furthermore, the radial degrees associated with
the curvature terms are always even in case of geodesic spheres and odd
in case of geodesic balls. Because of these differences, no proper
eliminations of terms appearing in $FineSIS_B(p)$ can be defined 
by the terms appearing in $FineSIS^{Ex}_S(p)$. Thus the proper 
spectral data considered for the unification of the ball-spectra 
with the sphere-spectra 
must be defined by the direct sum $FineSIS_S(p)\oplus FineSIS_B(p)$.

\begin{summary} {\bf (A) Case of spheres.}
On harmonic manifolds the Lichnerowicz identity $\mathcal Li(p)$
is not encoded into the space, $SIS^{ex}(p)$, of spectrally determined 
identities defined for the spectra of local geodesic spheres. 
Since this space, defined by the radial expansion coefficients of
functions $a_k(p,r)$ and $\int\theta^k(u_p,r)du_p$, 
is identical to the complete space of
spectrally determined identities, the Lichnerowicz identity
is not determined by the given spectral data, that is, it can not be
computed by the spectra of local geodesic spheres.
So is any identity $\mathcal F=\mathcal S+Q\mathcal Li$, 
where $\mathcal S\in SIS^{ex}(p)$ and $Q\not =0$. In \cite{am-s},
such spectrally undetermined identity is used to compute 
$||\nabla R||^2(p)||^2$.  

The greatest imperfection of $SIS(p)$ is that the radial expansion
coefficients $\int F_c(u)du_n$, whose integral-evaluation provides
the curvature expressions, are considered without any reference to
the radial degree $c$, thus this data can not be used for constructing
radial functions. (Such radial constructions can be implemented just by
functions $\int F_c(u)du_nr^c$. Also note that, in case of geodesic 
spheres, this degree is even.) As it is pointed out above, eliminations
using coefficients having different radial degrees result, on one hand
radial functions from which the considered term is not eliminated,
and, on the other hand, the curvature expression resulting by the formal
elimination does not appear as an expansion coefficient of a function 
being in $Span\{a_k(p,r)\}$ or 
$Span\{a_k(p,r),\int\theta^k(u_p,r)du_{np}\}$. In order to control this 
contradictory situation, refined spectral data        
$FineSIS(p)=\oplus_cSIS_c(p)$ resp. 
$FineSIS^{ex}(p)=\oplus_cSIS^{ex}_c(p)$, 
where $SIS_c(p)$ resp. $SIS^{ex}_c(p)$ consists of coefficients having
the same radial degree $c$ should be introduced to the computations. 
This data is defined after choosing the same
type of expansion for each function from the set $\{a_k(p,r)\}$ resp.
$\{\int\theta^k(u_p,r)du_{np}\}$.
The direct sum applied in these definitions
excludes the usage of coefficients of different radial degrees during
eliminations. Due to the identity 
$\int r^{2d}\theta^k(u_p,r)du_p=r^{2d}(\int\theta (u_p,r)du_{np})^k$ 
valid on harmonic manifolds, a larger, spectrally determined data 
$FineSIS^{ex}(p)\subset FineSIS^{Ex}(p)$ can be defined by the elements
of the larger function space  
\[
Span\{a_k(p,r),\int\theta^k(u_p,r)du_{np}\}\subset 
Span\{a_k(p,r),\int r^{2d}\theta^k(u_p,r)du_{np}\}.
\]
Although (\ref{denseq}) and the investigated $r^2$-coefficients of
$a_2(p,r)$ appear with different degrees, they can be brought together
to have the same radial degree by the above degree-adjusting operations.
Also note that the Lichnerowicz identity (which is not derived from
a radial function $l(p,r)$, thus no radial degree can be defined for it) 
can not be built into the space $FineSIS(p)$. By the above arguments, it
does not show up even in the rudimentary $SIS(p)$. Therefore, the ultimate
reason for $||\nabla R||^2(p)||^2$ is not determined by the spectra of
local geodesic spheres is that it can be computed just by the 
Lichnerowicz identity. 

{\bf (B) Case of balls.}
In the ball case, the AM\&S-paper proposes elimination of $L$ occurring
in the $r^3$-coefficient of $a_2(B_p(r))$ by the $L$ occurring in the 
$r^7$-coefficient of the spectrally determined function $vol(B_p(r))$.
Since these terms are defined for different radial degrees, by the
above explanation, this rudimentary elimination can not be applied
for these computations. An other problem is that the considered 
$r^3$-coefficient contains higher order powered terms which also must
be eliminated. For this operation one must use also 
the proper k-volumes, 
$Vol_k(p,r)=\int_{B_p(r)}\Theta^k_pdB_p(r)$, defined for $k\geq 1$.
Since these volumes are not determined by the standard volume
$vol(B_p(r))=Vol_0(p,r)$, they are  
not determined by the Dirichlet or Neumann spectra of balls either. 
Due to these facts, powered identities associated
with $\int C^3(u)du_n$ and $\int C(u)H(u)du_n$ become spectrally
undetermined. As a result, also $||\nabla R||^2$ becomes
spectrally undetermined.

Proper elimination of these terms can be implemented just by
using the adjusted volumes 
$Vol^{adj,\rho^{-4}}_k(p,r)=\int_{B_p(r)}\rho^{-4}\Theta^k_pdB_p(r)$,
where $k\geq 0$. But these adjustments further worsen the spectral
indeterminacy of the investigated terms. Even  
$Vol^{adj,\rho^{-4}}_0(p,r)$ becomes a spectrally undetermined quantity,
meaning, that the spectral determinacy of $L$ (regarding the spectra of
local geodesic balls) is a relative concept. It depends on the
degree of the radial expansion coefficient in which it shows up. In
the considered $r^3$-coefficient it is a spectrally undetermined term
and so are the other powered density-terms occurring there. This is
an other proof of the spectral indeterminacy of $||\nabla R||^2$.
These arguments also show that, in case of geodesic balls,  
extensions $FineSIS^{ex}(p)$ or $FineSIS^{Ex}(p)$ provide spectrally
undetermined quantities to the elements of $FineSIS(p)$. This is
the most drastic difference between the sphere- and ball-cases.

Quantity $||\nabla R||^2$ can be determined just by appropriate extended
data. Several of them are constructed above. 
All of them contain spectrally undetermined 
identities which play the very same role than what the Lichnerowicz 
identity
does when space $SIS^{ex}(p)$ of spectrally determined identities
is extend by the Lichnerowicz identity in order to determine 
$||\nabla R||^2$. It is also pointed out above that this quantity
is not determined even by any of the unified spectral data.   

\end{summary}

\section{Rudimentary uncertainty on Riemann manifolds.} 

The above considerations exhibit strong physical 
contents. Namely, one can explain by the given spectral data  
that how does the uncertainty principle manifest itself on
Riemannian manifolds. Here we consider only the rudimentary case. 
This idea works out on general Riemannian manifolds,
however, by using the available formulas, we start to explain it 
on harmonic manifolds. Despite this analogy, it should  
be pointed out that this theory is pure
geometric which can not be called physical, by any means. 
Although the constants and
and relations associated below with spectrally undetermined 
identities (such as the Lichnerowicz identity) strongly remind 
the Planck constant and the
corresponding Heisenberg relation, 
no physical interpretation for these 
identities exist at this point. This is why this theory is not associated
with any kind of physics at this point. 
It has no apparent relation even to the
exact physical models described in the last section.

The constants and relations mentioned above appear on the scene when
the distance of identities (such as the Lichnerowicz
identity) from $SIS(p)$ is computed. These computations are carried out
in a Hilbert space in the following way.
Consider all pre-functions $F(u)$ defined in terms of arbitrary
linear combinations of radial curvature
terms appearing in the AM\&S-expansions. In the expansions these terms
appear just in particular combinations but in generalized situations
all linear combinations are considered. The integrals of functions $F(u)$
define a linear space of equations endowed with the pre-Hilbert inner 
product
defined by integrating products $F_1(u)F_2(u)$ by $du$. By topological
closure, one defines an extended identity Hilbert-space $EIS(p)$. 
Now let $P(\mathcal L)$ be the spectrally determined 
identity obtained by projecting the Lichnerowicz identity onto  $SIS(p)$.
It can be obtained by integrating the projected function $P(L(u))$ by 
$du$. The identities in $SIS(p)$ are interpreted such that
they can be clearly heard by the given spectral data while the
other identities in $EIS(p)$ which are in the complement of $SIS(p)$ 
can be heard just together with noise. Such noisy identity 
is the Lichnerowicz identity, 
for which the noise wave is defined by $NW(L(u))=L(u)-P(L(u))$. 
The magnitude of noise is defined by the $L^2$-norm 
$||NW(L(u))||^2=LH(||\nabla R||^4,...,\hat R\mathring R,...)$ 
which is equal to an expression written in terms of the curvature 
terms appearing in the equations. Like in the physical Heisenberg 
relations, these terms are typically products of those appearing
in the primary equation
$\mathcal{NW}(L)$. Note that these quantities are defined by the
difference $L(u)-P(L(u))$, which makes the reminiscence with the
Heisenberg relations even stronger.
The Lichnerowicz-Planck constant, $h_{\mathcal L}$,
is defined by the minimum of values $||NW(L(u))||^2$ considered for all
$L(u)$ whose integration provides $\mathcal L$. Then the 
Lichnerowicz-Heisenberg relation is 
\begin{equation}
||NW(L(u))||^2=
LH(||\nabla R||^4, ||\nabla R||^2\hat R,||\nabla R||^2\hat R,
\hat R\mathring R,...)
\geq h_{\mathcal L},
\end{equation}
which clearly describes the positive distance of $\mathcal L$ from the
spectrally determined identities. 
 
A more general version of these relations can be defined such that
one replaces $P(L(u))$ by an arbitrary other function 
$F(u)\in Pre(SIS(p))$, by which the noise wave, $NW_F(L(u))$, is 
defined in the same way as by means of $P(L(u))$. 
Quantity $||NW_F(L(u))||^2$
measures the magnitude of noise when the Lichnerowicz identity is listened
by $F(u)$. Formally this function looks the same as the previous one, but
the noise is going to be bigger. That is, the left side of the general
HL-relation $||NW_F(L(u))||^2\geq h_{\mathcal L}$ usually gets
bigger in this situation.

These objects can be defined on arbitrary Riemann manifolds 
for any identity. The spectral data can be defined by local geodesic
spheres as well as balls. In order to avoid
the above discussed problems, the AM\&S-expansions must be 
established by the non-normalized density. Hilbert spaces 
$SIS(p)$ and $EIS(p)$ can be 
defined in the same way as they were on harmonic manifolds. 
It is also well known that Lichnerowicz established his
identity for general Riemannian manifolds, where it appears in a much
more complicated form, containing additional terms beyond those
appearing on Einstein manifolds. Anyhow, one can prove also in these 
most general cases that this general identity is not encoded into the
the spectra of local geodesic spheres or balls. More precisely we have:

\begin{summary} On a general Riemann manifold   
the Lichnerowicz identity, considered at a point $p$, is not determined
by the spectra of local geodesic spheres or balls concentrated at $p$. 
That is, this identity is not encoded into the spectral identity space  
$SIS(p)$ but it rather shows up as a spectrally undetermined
equation in $EIS(p)$. This fact can be explained also such 
that pre-functions defining the Lichnerowicz identity 
do not show up in $Pre(SIS(p))$ but in the larger pre-space
$ExPre(SIS(p))=Pre(EIS(p))$ where they have strictly positive 
distance from $Pre(SIS(p))$. The corresponding noise
waves, Plank constants, and Heisenberg relations can 
explicitly be computed also in these most general situations.

This mathematical uncertainty theory matches 
the physical one, by all means.
A curvature identity $\mathcal F$ defined by the integral $\int F(u)du_n$
of a pre-function is always a true identity. This logical value is 
independent from its relation to the spectral identity space $SIS(p)$.
The above quantities measure that in what extend can this well determined
identity be measured (or, recovered) by the given spectral data. It is
also clear that spectrally undetermined identities come from two sources.
Either by powering spectrally determined identities, or by considering such
identities like Lichnerowicz's which are nothing to do with pre-functions 
$F(u)\in Pre(EIS(p))$.
\end{summary}     

Several alternative versions of the above concepts
can be introduced as follows. For establishing the first one, 
consider a geodesic sphere or ball
of radius $R(p)$ about the center $p$. In the following definition the
data can either be the spectra of local geodesic spheres or balls 
concentrated at $p$. If the Taylor expansion 
$a_k(p,r)=\sum_{s=0}^\infty\frac 1{s!} A_{ks}(p)r^s$
of the local heat invariants is convergent for $r=R(p)$,
then, for any fixed $k$, heat invariant
$a_k(p,R)=\sum_{s=0}^\infty\frac 1{s!} A_{ks}(p)R^s$ (of the sphere
or ball of radius $R(p)$) can be considered as
a curvature identity, $\mathcal A_k(R(p))$, defined as an infinite sum 
(series) of identities $\frac {R^s}{s!}A_{ks}\in SIS(p)$ determined 
for the coefficients of this expansion. By this construction, one can
associate the identity space
$SIS(R(p))=\sum_{k=0}^\infty\mathcal A_k(R(p))$
to the spectrum of any fixed geodesic sphere or ball, which is just a very
thin subspace of the total space $SIS(p)$. For a field of geodesic 
spheres or balls, having radius $R(p)$ at a point $p$,
the $SIS(R(p))$ defines a Hilbert space bundle over the manifold, for 
which the Planck constants and
Heisenberg relations can be defined in the same way as
for $SIS(p)$. If an identity is not determined by $SIS(p)$, it is even
so regarding $SIS(R(p))$. This construction shows that the information
encoded into the spectral data is considerably weakened if, instead of
all local geodesic spheres or balls concentrated at $p$, it is provided
only by a single sphere or ball having center $p$.

Other concepts analogous to the AM\&S local heat invariants are the
so called {\it remote local heat invariants}, $\rho_k(p,\tau)$, defined
at a point $p$ of a compact Riemann manifold by means of the cut locus
$Cut(p)$. For a unit vector $u_p\in T_p(M)$, let $q(u_p)\in Cut(p)$ be
the point such that the geodesic starting out from $p$ into the direction
of $u_p$ intersects  $Cut(p)$ at $q(u_p)$. Parameterization, $\tau$,
on this geodesic, $c_{pq}(\tau )$, is chosen such that  
$c_{pq}(0)=p, c_{pq}(1)=q(u)$, furthermore, the speed vector 
$\dot c_{pq}(\tau )$ has constant length equal to the arc-length of  
$c_{pq}$. This parameter is called cut-locus-radius function. 
For a fixed $0\leq\tau\leq 1$, cut-locus-sphere
$CuS_p(\tau )$ of center $p$ and radius $\tau$ is defined by consisting of 
points which have
parameter $\tau$ on each of the above geodesics. Cut-locus-balls,
$CuB_p(\tau )$ are defined by points having parameters $\leq\tau$.
The remote local heat invariants are defined by the cut-locus-radial
expansion (i. e., $\tau$-expansion) of the heat invariants 
$a_k(CuS_p(\tau ))$ resp. $a_k(CuB_p(\tau ))$ defined on 
cut-locus-spheres resp. cut-locus-balls. Spectral identity spaces 
$SIS(Cut(p))$ and $SIS(Cut(\tau (p)))$, where $\tau$ denotes fixed 
cut-locus-radius, regarding cut-locus-spheres resp. balls can be
introduced in the same way as for geodesic spheres resp. balls.

The above constructions allow to introduce spectral identity spaces
also for compact sub-domains, $D$, with boundaries, $\partial D$, of
Riemannian manifolds $M^n$. For the sake of simplicity suppose that 
$\partial D$ is diffeomorphic to the Euclidean unit sphere $S^{n-1}$,
furthermore, let $p\in D$ be a point such that for any $q\in\partial D$ 
there exists a unique geodesic $c_{pq}\subset D$ joining the two points.
On such a geodesic, one can define the same parameterization, $\tau_p$, 
satisfying the very same properties described above. 
By replacing cut locus $Cut(p)$ with
boundary $\partial D$ in the above constructions, one can define
spectral identity space $SIS(p,\partial D)$ (where $\partial D$ 
appears as level set $\tau_p =1$) regarding the 
spectrum of Laplacian defined on $\partial D$ or any of the spectra
defined by particular boundary conditions on $D$. By continuous movement
of $p$, one can define this identity space together with the associated
Planck constants and Heisenberg relations at any point lying 
in the interior of $D$. Both the spectrally determined and undetermined
identities can precisely be described in terms of the Hilbert space 
bundles established in this construction.

\section{Global vs. local spectral investigations.}

On compact Riemann manifolds, $M$, the integrals of the 
infinitesimal heat invariants define the so called
averaged infinitesimal heat invariants. 
The exploration of relations between the global spectra and the averaged 
infinitesimal heat invariants seems to be an interesting question.
In this paper we establish just a single theorem
concerning the relation between the global spectra
and the volume of small geodesic balls and spheres. 
Next we prove that, contrary to the volume of the whole compact
manifolds, these local averaged volumes are not determined by the
spectrum of $M$.

By (\ref{boldA}), the fourth coefficient appearing in the
radial expansion of $\mathbf A_{c(r)}^k$ is a linear combination
of terms $TrR^2_{u_p}$ and  
$TrR^{\prime\prime}_{u_p}=Ricc(u_p,u_p)|u_p|u_p$, where $|u_p|u_p$
denotes two covariant differentiations regarding $u_p$.
Integration of $TrR^{\prime\prime}_{u_p}$ regarding $du_p$ at a fixed
points provides linear combinations of terms like $\Delta (Scal)(p)$ and
$(Ricc^{ab}_{|a})_{|b}$, whose integrals vanish on the whole compact
manifold, by the Stokes theorem. The same integral applied to 
the first term provides a constant time of $||Ricc||^2+(3/2)||R||^2$.

To establish the statement we consider  
also the one parametric families $(G,g_{\lambda (t)})$
of 0-isospectral manifolds constructed by Sch\"uth \cite{s} 
on compact Lie groups $G$ such as: 
\begin{eqnarray}
SO(m)\times T^2\,\, (m\geq 5),\, Spin(m)\times T^2\,\, (m\geq 5),
\,\,SU(m)\times T^2\, (m\geq 3),
\label{s1}
\\
SO(m)\, (m\geq 9),\,\, Spin(m)\, (m\geq 9),\,\, SU(m)\, (m\geq 6),
\,\, SO(8), \,\, Spin(8).
\label{s2}
\end{eqnarray}
She proved that distinct members
in a one-parametric family have distinct spectra on 1-forms by showing
that they have Ricci tensors of different norm. We use the
very same computational technique (she applied to prove this
latter statement) for establishing our statement.

Consider $a_2(G,g_{\lambda (t)})$ for a family, which, due to the 
isospectrality, is constant regarding $t$. Since $Scal$ is constant on 
$G$, furthermore, $vol(G)$ and $\int_GScal$ are spectrally determined,
also the first term, $\int_GScal^2$, in  $a_2(G,g_{\lambda (t)})$ is
constant regarding $t$. It follows then that the averaged volumes
of small geodesic balls and spheres can not be constant, otherwise
the above fourth coefficient in the expansion of the density together
with $a_2(G,g_{\lambda (t)})$ define a non-degenerated
system of linear equations which would determine constant values both for
$||Ricc||^2(t)$ and $||R||^2(t)$. This idea works out for all invariants
$Tr\mathbf A^k_{c(r)}$. By summing up we have

\begin{theorem}
The averaged volumes of small geodesic balls and spheres are generically
not determined by the spectra of compact Riemann manifolds.
This spectral indeterminacy is exhibited  
by the one parametric families $(G,g_{\lambda (t)})$
constructed by Sch\"uth \cite{s} on the compact Lie groups $G$ listed in
(\ref{s1}) and (\ref{s2}). The fourth coefficient in the power
series expansion of the volumes of small geodesic balls resp. spheres 
depends on parameter $t$. This
coefficient is expressed in terms of $||Ricc||^2$ and other lower order
curvature invariants. Sch\"uth originally proved that quantity 
$||Ricc||^2$ is different on the distinct members of a family, therefore,
she concluded, they are not isospectral on 1-forms.
\end{theorem}

This statement can be established, by the same proof, for 
$Tr\mathbf A^k$-volumes of geodesic balls and spheres where the 
integration can be defined either with $\theta du$ or just with $du$.
This statement is highly expectable also on the Gordon-Wilson 
examples \cite{gw}
of isospectral manifolds having different local geometries as well 
as on those
derived from these examples. Since these manifolds have non-constant
scalar curvature in general, this version of the theorem can not be
established by using only the fourth coefficient of the above power series 
expansion.

\section{Isospectralities and physical symmetries.}

Both my isospectrality examples and the Gordon-Wilson
examples \cite{gw} (together with those derived from the GW-constructions,
see a survey on them in \cite{s}) 
arise from 2-step nilpotent Lie groups. Yet, in terms of spectral 
stability of small geodesic balls and spheres 
these two types of examples are the polar opposites of each other.
One of the idiosyncrasies of the GW-examples and their relatives is
the wide range of changing invariants during the continuous isospectral
deformations introduced in their constructions. Whereas, 
these quantities are not changing 
during the discrete isospectral deformations introduced in my 
constructions.
In my isospectrality families the members share even the same volume
function $\theta (r)$, therefore, also the volumes of geodesics spheres 
resp. balls having the same radius must be identical. Thus the above
theorem is a clear demonstration of these arguments. Later, after 
introducing the technical definitions, it is more clearly pointed out
that what causes these significant differences between these 
two different type of constructions.
 
The isospectrality exhibited by my examples has a deep physical meaning. 
Namely, it is equivalent to the C-symmetry known in quantum theory.
This physical connection is completely established in \cite{sz8}.
Actually, the bulk of this whole section is going to be an outline 
of some of the results developed in \cite{sz8}. Readers interested
in details may consult with this paper, however, the following review tends
to be self-contained as much as it is possible.

Before going into technical details, we briefly 
describe this physical connection by non-technical words.  
As it is pointed out below, on the Riemann manifolds used for these 
constructions the Laplacian is nothing but the quantum Hamilton 
operator of particle systems described in elementary
particle physics. The action of the intertwining operator constructed 
for establishing the isospectrality can be interpreted such that it 
exchanges some of the particles for their anti 
particles. Thus the mathematically established isospectrality
can be paraphrased such that the spectrum of an elementary
particle system, that is, the possible energy levels on which 
the system can
exist, is not changing if some of the particles are exchanged for 
their anti particles. 

This paraphrase is a clear manifestation of the C-symmetry principle
introduced in physics in the following general form: 
``The laws are the same for particles and antiparticles."
It should be emphasized, however, that this general 
principle will be pointed out only regarding the spectrum of the Hamilton
operator. 
In other words, the isospectralities exhibited by my examples
are clear manifestation of the {\it spectral C-symmetry principle.}

A preliminary description of the endomorphisms exchanging the 
particles for their antiparticles is as follows.
The two-step nilpotent Lie groups used for these constructions are
defined on the (X,Z)-space, 
$\mathbb R^k\times\mathbb R^k=\mathcal X\times\mathcal Z$ by a linear 
space, 
$E_{skew}(\mathcal X)$, of skew endomorphisms
acting on the X-space. The members of an isospectrality family
are defined on the same (X,Z)-space 
by a corresponding family of endomorphism spaces where for any
two members, $E_{skew}(\mathcal X)$ and
$E^\prime_{skew}(\mathcal X)$, there exists an involutive 
orthogonal transformation $\sigma : \mathcal X\to\mathcal X$ commuting
with all these endomorphisms and the endomorphisms belonging to
$E^\prime_{skew}(\mathcal X)$ can be obtained by composing the elements of 
$E_{skew}(\mathcal X)$ by $\sigma$, that is,
$E^\prime_{skew}(\mathcal X)=\sigma\circ E_{skew}(\mathcal X)$.
These conditions imply the existence of an orthogonal direct sum 
$\mathcal X=\mathcal X^{(a)}\oplus\mathcal X^{(b)}$ where both
subspaces are invariant under the actions of all endomorphisms 
and the first component is fixed under the action of $\sigma$,
while it is $-id_{\mathcal X^{(b)}}$ on the orthogonal complement 
$\mathcal X^{(b)}$. 
It is explained later that these endomorphisms 
serve as angular momenta for particles orbiting in $\mathcal X^{(a)}$
resp. $\mathcal X^{(b)}$. These endomorphisms reverse the angular
momentum for particles living on $\mathcal X^{(b)}$ and preserve it 
for those orbiting in $\mathcal X^{(a)}$. The operators intertwining the
Laplacians are defined by means of these exchange-endomorphisms.  
Let it be mentioned yet that these constructions 
allow just discrete isospectral deformations.

By using the same two-step nilpotent Lie group in two different ways,
these examples are constructed on two different type of manifolds. 
The first ones are torus bundles defined over the 
X-space by factoring the center, $\mathcal Z$, 
by a Z-lattice $\Gamma_Z$, and the others are ball resp.
sphere bundles over the X-space defined by considering appropriate
Z-balls resp. Z-spheres in $\mathcal Z$. For the sake of simplicity,
next we describe only the ball-bundle cases. The spectral investigations
must be established on the function spaces defined on these manifolds. 
Actually both function spaces can be considered on the same non-factorized
manifold. Namely, in the first case, it is the space of 
$\Gamma_Z$-periodic functions, while in the second case it is the function
space satisfying an arbitrarily prescribed boundary condition 
on the boundary manifold which is a Z-sphere bundle. Both function space
can be written up by appropriate Z-Fourier transforms. In the first
case it is the discrete Z-Fourier transform defined by the Z-lattice.
The big advantage of this transform is that it
separates the X- and Z-variables from each other. The other 
Z-Fourier transform introduced for investigating Z-ball bundles is the so 
called twisted Z-Fourier transform. This name indicates that this
is a much more
complicated version of the Z-Fourier transforms where the X- and 
Z-variables can not be separated from each other. This transform 
involves mixed, so called twisting functions which establish the 
fundamental bonds between the X- and Z-spaces.

One of the most important observation regarding these two 
exact mathematical models is that with the help of them 
the electromagnetic, the
weak-, and the strong-nuclear forces can be established within a 
unified framework. The main unifying principle is that
these forces can be described by the eigenfunctions of the very same
Laplacian such that the distinct forces emerge on 
distinct invariant 
subspaces of this common quantum operator. This very same Laplacian
is the Laplacian acting on the Lie group and the two distinct function
spaces are the $\Gamma_Z$-periodic function space resp. those defined
for Z-ball bundles. After its action on the Fourier integral formula,
the very same Laplacian appears in drastically different ways behind
the integral sign. In the first case, it turns into a 
Ginsburg-Landau-Zeeman operator of charged particles orbiting
in complex planes in magnetic fields perpendicular to the planes, while,
in the second case, it appears as a sum of a scalar operator
which can be identified as the quantum Hamilton operator
of electroweak interactions and a new type of spin operator which
can be identified as the quantum Hamilton operator of strong force
interaction keeping the nucleus together. The particles attached 
to the two function spaces are distinguished by calling
them particles having no interior resp. those having interior.  

The Hamilton operators are defined in this theory on the nilpotent
groups.
The corresponding wave and Schr\"odinger operators 
emerge in the Laplacians of the static resp. 
solvable extensions of these nilpotent groups.
The latter manifolds are endowed with a natural 
invariant indefinite metric of 
Lorentz signature. Thus, these new exact 
mathematical models 
provide a relativistic theory for elementary particles. The
above discussed functions defined by Z-Fourier transforms
appear in the explicit solutions of these wave operators.
These wave functions strongly remind those introduced by de Broglie 
in classical wave mechanics (cf. \cite{p}, volume 5). Actually, they are 
the only appropriate adaptations
of the original de Broglie wave functions to the new mathematical
models. They carry over also the original de Broglie theory to understand 
the much more complicated physical situation inherent 
in the new models. 
This theory establishes infinitely many
non-equivalent models for which even classification is possible. The
particle systems attached to them behave exactly like
those introduced by the familiar standard model of 
elementary particle physics.

\subsection{Technicalities on 2-step nilpotent Lie groups.}

A 2-step nilpotent metric Lie algebra, 
$\{\mathcal N ,\langle ,\rangle\}$, is defined on a real vector space 
endowed with a positive definite inner product.
The name indicates that the 
center, $\mathcal Z$, can be reached
by a single application of the Lie bracket, thus its second application
always results zero. 
The orthogonal complement of the center is denoted by $\mathcal X$. 
The Lie bracket operates among these subspaces according to
the following formulas: 
\begin{equation}
[\mathcal N,\mathcal N]=\mathcal Z\quad ,\quad 
[\mathcal N,\mathcal Z]=0\quad ,\quad 
\mathcal N=\mathcal X\oplus\mathcal Z=\mathbb R^{k}\times\mathbb R^l.
\end{equation}
Spaces $\mathcal Z$ and $\mathcal X$ are called Z- and
X-space, respectively. 

Up to isomorphisms, such a Lie algebra is uniquely
determined by the linear space, $J_{\mathcal Z}$, of skew endomorphisms 
$J_Z:\mathcal X\to\mathcal X$ defined for Z-vectors
$Z\in\mathcal Z$ by the formula 
\begin{equation}
\label{brack}
\langle [X,Y],Z\rangle =\langle J_Z(X),Y\rangle , 
\forall Z\in\mathcal Z.
\end{equation}
This statement means that for an orthogonal direct sum, 
$\mathcal N=\mathcal X\oplus\mathcal Z=\mathbb R^{k}\times\mathbb R^l$,
a non-degenerated linear map, 
$\mathbb J:\mathcal Z\to E_{skew}(\mathcal X)\, ,\, Z\to J_Z$,
from the Z-space into the space of skew endomorphisms acting on the 
X-space, defines a 2-step nilpotent metric Lie algebra on 
$\mathcal N$ by 
(\ref{brack}). Furthermore, an other non-degenerated linear map 
$\tilde{\mathbb J}$ having the same range
$\tilde J_{\mathcal Z}=J_{\mathcal Z}$ as 
$\mathbb J$ defines a Lie algebra which is isomorphic to the previous one.
If isometric isomorphism is required, then 
$\tilde{\mathbb J}^{-1}\circ\mathbb J$ must be an orthogonal 
transformation on the Z-space.

With the help of these technical definitions one can more clearly 
explain why do the GW-type constructions produce completely 
different examples from those I
constructed by the above exchange-endomorphism $\sigma$. 
In the GW-case the isospectral deformation is
defined such that, for any fixed $Z$, the spectrum of the bilinear
map $\langle J_Z(X_1),J_Z(X_2)\rangle$ (which defines a unique linear
map on the Euclidean X-space) is not changing, but, the  spectrum of 
$\langle J_{Z_1}(X),J_{Z_2}(X)\rangle$ (defined, now, on the
Euclidean Z-space for an arbitrary fixed $X$) is changing. 
Other characteristic
features are that these maps have different eigenvalues, furthermore, the
Z-space has dimension $2$. The changing spectrum of  
$\langle J_{Z_1}(X),J_{Z_2}(X)\rangle$ gives rise to the wide range
of changing invariants (including also spectral invariants) defined 
for Z-balls, Z-spheres, geodesic balls, and geodesic spheres. This is the
ultimate reason why the GW isospectrality examples are established just on
Z-torus bundles but not on Z-ball and Z-sphere bundles, where they are 
actually non-isospectral. 

Whereas, the spectra both of
$\langle J_Z(X_1),J_Z(X_2)\rangle$ and  
$\langle J_{Z_1}(X),J_{Z_2}(X)\rangle$ are not changing during the
discrete deformations I introduced for the constructions. This is why
my isospectrality examples live both on Z-torus bundles and Z-ball resp.
Z-sphere bundles. Due to these differences, the following formulas
established in my case resist any attempt to define them on the
GW-examples. In this paper we describe only such cases where the spectra
of the above maps consist only the same values. Thus difficulty arises
already in the first step when one tries to extend them to the
GW-type constructions, where these spectra consist of different 
eigenvalues. But the ultimate reason for this extension breaks down 
is that the spectrum of $\langle J_{Z_1}(X),J_{Z_2}(X)\rangle$ is changing
during the continuous deformations they introduced for 
their constructions. These difficulties strongly indicate that the 
distinguishing characteristics exhibited by these two different cases 
can be  explored by mutually distinct methods.  
This field can not be explored by a single method with the help
of which one would be able to control both the changing and 
the non-changing invariants exhibited by these two different type of 
constructions.

Important particular 2-step nilpotent Lie groups are the
Heisenberg-type groups \cite{k}
defined by endomorphism spaces $J_{\mathcal Z}$
satisfying the Clifford condition $J^2_Z=-\mathbf z^2id$, where 
$\mathbf z=|Z|$ denotes the length of Z-vectors. 
These groups are attached to Clifford modules 
(representations of Clifford algebras). 
The well known classification 
of these modules provides  
classification also for the Heisenberg-type groups. According to this
classification, the X-space of a H-type group is an $(a+b)$-times
Cartesian product of a smaller space 
$\mathcal Y=\mathbb R^{n_l}$, which is 
endowed with an $l$-dimensional endomorphism space $j_{\mathcal Z}$
such that the endomorphisms acting on 
$\mathcal X=\mathcal Y^a+\mathcal Y^b$ are defined by
$J_Z=j_Z\times\dots\times j_Z\times -j_Z\times\dots\times -j_Z$. Dimension
$n_l$ depends only on $l$. 

By means of the exponential map, also
the group can be considered to be defined on $\mathcal N$.
That is, a point is denoted by $(X,Z)$ also on 
the group.   
Metric tensor, $g$, is defined by the left invariant extension of 
$\langle ,\rangle$ onto the group $N$.

Although most of the results described below extend to general 
2-step metric nilpotent Lie groups, next only H-type groups
will be considered. 
On these groups, the Laplacian appears as follows:
\begin{equation}
\label{Delta}
\Delta=\Delta_X+(1+\frac 1{4}|X|^2)\Delta_Z
+\sum_{\alpha =1}^r\partial_\alpha D_\alpha \bullet,
\end{equation}
where $\{e_\alpha\}$ is an orthonormal basis on the center (Z-space)
and $D_\alpha\bullet$ denotes directional derivatives along
the vector fields 
$X\to J_\alpha (X)=J_{e_\alpha}(X)$, furthermore, $\mathbf x=|X|$
denotes the length of X-vectors.

In the next sections, it is described in technical terms that how this  
operator manifests itself as the Hamilton operator of different 
elementary particle systems. A non-technical preview of this
physical interpretation is as follows.
As it is indicated above, the main idea is that these
systems are attached to various invariant function subspaces of $\Delta$
which can be divided into two major classes. 
The first class is defined by Z-torus bundles 
(alias Z-crystals), where the corresponding function space consists of 
functions which are periodic regarding the Z-lattice $\Gamma_Z$ defining
the Z-torus bundle. 
The particles  attached to such a
function spaces are considered to be point particles having no
interior. They can show up at the lattice points of $\Gamma_Z$, where
the angular momentum is defined for them by means of $J_Z$. 
As it will turn out, the $\Delta$ appears on
this function space as the Ginsburg-Landau-Zeeman operator of
charged particles. Thus the natural interpretation for these systems
is that they consist of electrons, positrons, 
and electron-positron-neutrinos. Contrary to these cases, 
the other type of particles,
emerging on function spaces defined on Z-ball bundles by
Dirichlet or Neumann boundary conditions, do have interior which space is
represented by the interior of the Z-balls. The X-space represents
always the exterior word. The Laplacian appears on these Z-ball-bundles
as the sum of Hamilton operators of electro-weak and
strong-force interactions. A technical description of these enormous
differences between Z-torus and Z-ball bundles is as follows.

\subsection{Z-crystals modelling Ginsburg-Landau-Zeeman operators.}
The Z-torus bundles are defined by factorizing,  
$\Gamma\backslash H$, a two-step nilpotent group $H$ by a Z-lattice,
$\Gamma=\{Z_\gamma\}$. The name indicates that this lattice is defined  
only on $\mathcal Z$ and not on the whole $(X,Z)$-space. 
Such a factorization defines a Z-torus bundle over the X-space.
The natural Z-Fourier decomposition,
$L^2_{\mathbb C}:=\sum_\gamma W_\gamma$, of the $L^2$ function space 
is defined such that subspace
$W_\gamma$ is spanned by functions of the form
\begin{equation}
\label{diskFour}
\Psi_\gamma (X,Z)=\psi (X)e^{2\pi\mathbf i\langle Z_\gamma ,Z\rangle}.
\end{equation}
Note that each $W_\gamma$ is invariant under the 
action of $\Delta$, more precisely, we have:
\begin{eqnarray}
\Delta \Psi_{\gamma }(X,Z)=(\lhd_{\gamma}\psi )(X)
e^{2\pi\mathbf i\langle Z_\gamma ,Z\rangle},\quad
{\rm where}
\\
\lhd_\gamma
=\Delta_X + 2\pi\mathbf i D_{\gamma }\bullet 
-4\pi^2|Z_\gamma |^2(1 + \frac 1 {4} |X|^2).
\end{eqnarray}
In terms of parameter $\mu =\pi |Z|_\gamma $, 
this operator can be written
in the form
$
\lhd_{\mu}=
\Delta_X +2 \mathbf i D_{\mu }\bullet -\mu^2|X|^2-4\mu^2.
$
Although it is defined in terms of the X-variable,
this operator is not a sub-Laplacian which is resulted by a submersion. 
It should be considered as restriction of the total Laplacian onto 
the invariant subspace $W_\gamma$.
Actually, the Z-space is represented by the constant $\mu$
and operator $D_\mu\bullet$. A characteristic feature of this 
restricted operator is that it involves only a single endomorphism,
$J_{Z_\gamma}$.

On a Z-crystal, $B_R\times T^l$, for given
boundary condition $Af^\prime (R^2)+Bf(R^2)=0$,
the eigenfunctions
of $\lhd_{\mu }$ can be represented in terms of $n^{th}$-order 
complex valued spherical harmonics
$\mathtt H^{(n,m)} (X)$
in the form $f(\langle X,X\rangle )\mathtt H^{(n,m)} (X)$, where the radial
function $f$ is an eigenfunction of the ordinary differential operator
\begin{equation}
\label{Lf_lambda}
(\large{\Diamond}_{\mu,\tilde t}f)(\tilde t)=
4\tilde tf^{\prime\prime}(\tilde t)
+(2k+4n)f^\prime (\tilde t)
-(2m\mu +4\mu^2(1 +
{1\over 4}\tilde t))f(\tilde t).
\end{equation}
Degree $m$ is defined such that 
$\mathtt H^{(n,m)} (X)$ is simultaneously  eigenfunction also 
of $\mathbf iD_\mu\bullet$ with eigenvalue $m\mu$. 

Polynomials $\mathtt H^{(n,m)} (X)$ can be constructed as follows. 
Consider a complex basis $\mathbf B=\{Q_1,\dots ,Q_{k/2}\}$ regarding
the complex structure $J_{Z_{\mu u}}$ and an $n^{th}$-order polynomial
$\prod z_i^{p_i}(Z_{\mu u},X)\overline z_i^{q_i}(Z_{\mu u},X)$, where 
$z_i(Z_{\mu u},X)=\langle Q_i+\mathbf iJ_{Z_{\mu u}}(Q_i),X\rangle$ and
$\sum (p_i+q_i)=n$. Since  $z_i(Z_{\mu u},X)$ is an eigenfunction of
$\mathbf iD_\mu\bullet$ with eigenvalue $-\mu$, the whole
polynomial is also an eigenfunction with eigenvalue
$m\mu$, where $\sum (-p_i+q_i)=m$. However, this homogeneous polynomials
are not harmonic regarding $\Delta_X$ and their restrictions do not define
spherical harmonics on the unit X-sphere. In order to have the 
spherical harmonics, these spherical but non-harmonic functions should
be projected into the space of $n^{th}$-order spherical harmonics.
These projections, $\Pi_X^{(n)}$, are extensively described in my papers
\cite{sz4}-\cite{sz8} (see for instance Section 6.4 in \cite{sz8}). 
They can be represented as uniquely
determined polynomials of the spherical Laplacian $\Delta_S$. Since
operator $\mathbf iD_\mu\bullet$ commutes with this Laplacian,
projection  
$\Pi_X^{(n)}((\prod z_i^{p_i}\overline z_i^{q_i}))(Z_{\mu u},X)$
really provides a desired polynomial $\mathtt H^{(n,m)} (X)$.
Since these projections are onto whose kernels are formed exactly 
by the lower order polynomials, all these polynomials can be obtained
by this construction.

In the 2D-case, operator $\lhd_{\mu}$ can be transformed 
to the the well known Ginsburg-Landau-Zeeman 
operator \cite{bo,ll,p} 
\begin{equation}
\label{land}
-{\hbar^2\over 2m}\Delta_{(x,y)} -
{\hbar eB\over 2m c\mathbf i}
D_z\bullet
+{e^2B^2\over 8m c^2}(x^2+y^2)
\end{equation}
of a charged particle orbiting in the $(x,y)$-plane in constant
magnetic field directed toward the z-axis
by choosing $\mu ={eB/2\hbar c}$ and  
multiplying the whole operator with 
$-{\hbar^2/2m}$.
For k-dimensional X-spaces,
number $\kappa =k/2$ means the {\it number of particles}, and, 
endomorphisms $j_Z$ and $-j_Z$ in the above formulas are attached to 
systems electrons resp. positrons. More precisely, by the classification 
of H-type groups, these endomorphisms are acting on the irreducible
subspaces $\mathbb R^{n_l}$ and the system is interpreted such that 
there are $n_l/2$ particles of the same charge orbiting on complex planes
determined by the complex structures $j_{Z_u}$ resp. $-j_{Z_u}$
in constant magnetic fields. 

Note that this operator contains
also an extra constant term $4\mu^2$, which, after establishing
the corresponding wave operators, can clearly be explained 
as the total energy of neutrinos accompanying 
the electron-positron system \cite{sz8}. 
(This energy term is neglected in the original
Ginsburg-Landau-Zeeman Hamiltonian.) 
Thus these mathematical models 
really represent systems of charged particles and antiparticles accompanied
with electron-positron-neutrinos known in elementary particle physics
\cite{v,w1,w2}.
Let it be mentioned yet that operator $D_\mu\bullet$ is associated with 
the magnetic dipole resp. angular momentum operators of classical quantum 
theory. They are defined for the lattice points separately. Thus, while
wandering on the lattice points, these point-like particles receive
their angular momenta at the lattice points where they stay on.

Finally, let the isospectrality question be clarified. An isospectrality
family is defined by Heisenberg type groups $H_l^{(a,b)}$ having the
same $l$ and $(a+b)$, that is, they share the same X-space
$\mathcal X=\mathcal Y^a+\mathcal Y^b$ and Z-space $\mathcal Z$.
They differ from each other just by the decomposition of the X-space
and the action of $J_Z$ on these components. When 
$H_l^{(a+b,0)}$ is compared with $H_l^{(a,b)}$, then the exchange
endomorphism $\sigma$ described above is defined by $id_{/\mathcal Y^a}$
resp. $-id_{/\mathcal Y^b}$ on the components of the above decomposition.
On $\mathcal X{(b)}=\mathcal Y^b$, this exchange endomorphism defines 
the angular momentum
of $H_l^{(a,b)}$ by the negative of the angular momentum defined for
$H_l^{(a+b,0)}$. One can interpret this as switching the sign of charge
of the particles on this component. Now consider an X-ball around
the origin and restrict the torus bundle onto this ball. For both bundles,
the normal
vector at a boundary point is the the radial unit X-vector. For both
metrics consider the same complex basis $\mathbf B$ and the same
polynomials in terms of the complex structures defined for these
two metrics. 
Since both define the same radial Laplacian with the same boundary
conditions, these two metrics on the considered sub-bundles must be 
isospectral regarding any of the boundary conditions. This isospectrality
can be seen also by observing that, for any fixed $Z_\gamma$, there exist
an orthogonal transformation on the X-space which conjugates 
$J^{(a+b,0)}_{Z_\gamma}$ to $J^{(a,b)}_{Z_\gamma}$, therefore, it
intertwines the Laplacians of the two metrics term by term 
along with the boundary conditions.

If basis $\mathbf B_\gamma$ is chosen by picking the vectors always
from subspaces $\mathcal Y=\mathbb R^{n(l)}$, this isospectrality
is a clear manifestation of the spectral C-symmetry. For a general
basis, however, it exhibits also internal symmetries. 
In elementary particle physics, this name
was chosen to indicate that one think about internal symmetries as having
to do with the intrinsic nature of the particles, rather than their
position or motion. You can think of each particle as carrying a 
little dial, with pointer that points in directions marked ``electron"
or ``neutrino" or ``photon" or ``W" or anywhere in between. 
The internal symmetry says that the laws of nature take the same form
if we rotate the markings on these dials in certain ways. If a basis 
vector $Q_i=Q^{(a)}_i+Q^{(b)}_i$ is lying neither in $\mathcal Y^{(a)}$
nor in $\mathcal Y^{(a)}$, the exchange endomorphism effects only
$Q^{(b)}_i$, that is, the particle antiparticle exchange is just partial
and not complete. By this explanation, this isospectrality is a clear 
manifestation of the spectral C/I-symmetry.  

\subsection{Extended particles occupying Z-ball bundles.} 
Contrary to the
above Z-crystal cases, the Laplacian on Z-ball
and Z-sphere bundles can be identified with the Hamilton operators 
of particles to which interior can be attributed. In elementary 
particle physics
such particles occur in the nucleus where the constituents are held 
together by the electroweak and strong forces. These forces are exhibited
by the idiosyncratic appearance of the Laplacian on these bundles.
Namely, it decomposes into a scalar operator and a 
non-standard spin operator (called also {\it roulette operator}) where
the scalar operator represents the electroweak force while the 
roulette operator corresponds to the strong force.
     
The ball$\times$ball- and ball$\times$\-sphe\-re-type 
manifolds used to these investigations emerged first in the spectral 
constructions performed in \cite{sz4}-\cite{sz6}. These manifolds 
are defined by appropriate smooth fields of Z-balls resp. Z-spheres 
of radius $R_Z(|X|)$ 
over the points of a fixed X-ball $B_{\mathcal X}$ 
of radius $R_{\mathcal X}$. 
Note that radius $R_Z(|X|)$ depends just on the
length, $|X|$, of vectors $X\in B_{\mathcal X}$ over which the 
Z-balls resp. Z-spheres are considered. The boundaries
of these manifolds are the so called
sphere$\times$ball- resp. sphere$\times$sphere-type manifolds.
Comparing with the Z-crystals, the difference between the 
two type of bundles is 
that one considers Z-balls resp. Z-spheres instead of the 
Z-tori used in the previous construction. In the
isospectrality investigations the compact domains
corresponding to $R_{\mathcal X} <\infty$ play the 
primary interest.
In physics, however, the non-compact bundles corresponding to 
$R_{\mathcal X} =\infty$ (that is, which are defined 
over the whole X-space) become the most important cases. In the following
considerations both cases will be investigated. The details will be
provided in this paper just for Z-ball bundles and not for 
Z-spheres bundles. 

The main difference between the Z-crystals and Z-ball bundles is that 
the computations in the latter case can not be reduced to a single 
endomorphism but they must be established for the complete
operator $\mathbf M=\sum\partial_\alpha D_\alpha\bullet$. This operator 
includes the angular momentum endomorphisms $J_Z$ with respect 
to any Z-direction. This complication gives rise to a
much more complex mathematical and physical situation
where both the exterior and the interior life of particle systems 
exhibit them-self on a full scale.

The above argument implies that this mathematical situation 
can not be described by the discrete Z-Fourier
transform applied on Z-crystals. 
In this case one considers
a fixed complex basis $\mathbf B$ together with the complex coordinate 
systems $\{ z_i(K_u,X)\}$ defined, on the X-space, by the complex
structures $J_{K_u}$. Then, 
the Z-Fourier transform is defined on the whole center, 
$\mathcal Z=\mathbb R^l$, by
\begin{equation}
\label{newDiracwave}
\int_{\mathbb R^l}A(|X|,K)\Pi_X^{(n)}(\prod z_i^{p_i}(K_u,X)
\overline z_i^{q_i}(K_u,X))e^{\mathbf i\langle K,Z\rangle}dK.
\end{equation}
It should be pointed out that a fixed $\mathbf B$ can serve as
a complex basis 
only for almost every $K_u$, which form an everywhere 
dense open subset on the unit Z-sphere. However, the polynomials
are well defined analytic functions even at those $K_u$'s with respect
to which the $\mathbf B$ is not a complex basis. 
This formula is well defined if for any fixed $|X|$ function
$A(|X|,K)$ is of class $L^2$ regarding the $K$-variable.

This so called
twisted Z-Fourier transform does not separate but rather binds 
the Z-space and the X-space together. This strong bond is established
by the polynomials $\prod z_i^{p_i}\overline z_i^{q_i}$ which depend
both on the X- and K-variables. 
It is proved in Section 6.2 of \cite{sz8} that the complex 
valued functions  
considered behind the integral sign 
for all possible powers satisfying $\sum_i(p_i+q_i)=n$
span an everywhere 
dense subspace, $\mathbb Tw_{\mathbf B}^{(n)}$, of the straightly defined 
complete function space, $\mathbb{S}t^{(n)}$, 
which can be introduced by extending
basis $\mathbf B$ into a basis 
$\tilde{\mathbf B}=\{\tilde Q_1,\dots ,\tilde Q_k\}$
of the whole X-space and replacing the above complex polynomials
by the real polynomials 
$\prod\langle \tilde Q_i,X\rangle^{\alpha_i}$, where $\sum{\alpha_i}=n$.
While $\mathbb {T}w_{\mathbf B}^{(n)}$ depends on $\mathbf B$, 
this straightly defined function space
is independent from the choice of the basis  $\tilde{\mathbf B}$. 
This function space naturally emerges on the Cartesian product 
$\mathcal X\times\mathcal Z$ as Cartesian product of the corresponding
function spaces. 
The same statements are true for the Z-Fourier transforms, 
$\mathbb F_Z(\mathbb {T}w_{\mathbf B}^{(n)})$ and
$\mathbb F_Z(\mathbb {S}t^{(n)})$,
of these function spaces. When these function spaces are defined
for particular powers $p_i$ and $q_i$, they are denoted by
$\mathbb F_Z(\mathbb {T}w_{\mathbf B}^{(p_i,q_i)})$ resp.
 $\mathbb F_Z(\mathbb {S}t^{(p_i,q_i)})$.

It is also observed in [Sz8] (cf. Theorem 6.1) that, 
for a function represented
by twisted Z-Fourier transform, $A(|X|,K)$ is a uniquely
determined 
$L_K^2$-function. Thus, by fixing a basis $\mathbf B$ and 
two H-type groups $H^{(a,b)}_l$ and 
$H^{(a^\prime ,b^\prime)}_l$ satisfying $(a+b)=(a^\prime +b^\prime)$,
there exists a well defined one to one correspondence 
$\kappa_{\mathbf B}:\mathbb F_Z(\mathbb {T}w_{\mathbf B}^{(n)})
\to \mathbb F_Z(\mathbb {T^\prime}w_{\mathbf B}^{(n)})$
that maps a
function expressed by $A(|X|,K)$ and complex structures 
$J_{K_u}^{(a,b)}$ to functions where just the complex structures are 
exchanged for $J_{K_u}^{(a^\prime ,b^\prime )}$. 
Let it be pointed out again that this map
operates on functions defined by the twisted Z-Fourier transform and not
on the ones standing behind the integral sign. This map has a
unique extension defining a bijection between 
the straightly defined function spaces.  

In order to check out if this map is an intertwining operator,
let the Laplacian be act on the twisted Z-Fourier transform formula.
Like for Z-crystals, operator 
$\mathbf M=\sum\partial_\alpha D_\alpha\bullet$
appears behind the integral sign in the form
$\mathbf i|K|D_{K_u}\bullet$, which acts only on the polynomial part,
resulting
$|K|m$, where $m=\sum (-p_i+q_i)$. Term involving $\Delta_Z$ appears
inside as $4\pi^2|K| ^2(1 + \frac 1 {4} |X|^2)$. Finally, the radial
Laplacian in $\Delta_X$ acts on X-radial functions by radial
differentiations and multiplications with radial functions,
furthermore, the action of the X-spherical Laplacian is nothing
but multiplication with the corresponding eigenvalue. Thus the
above operator really intertwines the Laplacians (cf. these details in
Section 8 of \cite{sz8}).   
  
Although the considered functions are everywhere dense in the 
straightly defined complete function space, the 
boundary conditions can not be computed by them. The main problem is
that, regarding the K-variable, function  
$A(|X|,K)\Pi_X^{(n)}(\prod z_i^{p_i}(K_u,X)\overline z_i^{q_i}(K_u,X))$
is not the multiple of a single $s^{th}$-order spherical harmonic by
a K-radial function even if one uses functions 
$A(|X|,K)=\phi(|X|,|K|)\varphi (K_u)$, where
$\varphi (K_u)$ is a spherical harmonic of degree $f=deg(\varphi )$.
This claim follows from the fact that function 
$\prod z_i^{p_i}(K_u,X)\overline z_i^{q_i}(K_u,X)$ is not derived
from a homogeneous polynomial but, after performing the indicated 
powering imposed on functions
 $z_i(K_{u},X)=\langle Q_i+\mathbf iJ_{K_{u}}(Q_i),X\rangle$ and
its conjugate, it appears as a proper sum of the form:
\begin{eqnarray}
\prod z_i^{p_i}(K_u,X)\overline z_i^{q_i}(K_u,X)=
\\
\sum_{(a_1\dots a_n)}\prod_i\langle Q_i,X\rangle^{n_i-a_i}
\langle\mathbf iJ_{K_{u}}(Q_i),X\rangle^{a_i}
=\sum_{a=\sum a_i=0}^nR^{(a)}(K_{u},X),
\nonumber
\end{eqnarray}
where $n_i=(p_i+q_i)$, $n=\sum n_i$, and index $a=\sum a_i$ indicates that 
how many complex structures 
$J_{K_{u}}$ appear in the term determined by
the exponents $(a_1\dots a_n)$. Sub-sum  
$R^{(a)}(K_{u},X)$ consists exactly those terms where this number is $a$. 
Although
these functions are derived from $a$-homogeneous functions,
neither they nor $\varphi (K_u)R^{(a)}(K_{u},X)$ are 
spherical harmonics regarding $K_{u}$. One obtains the
desired $s^{th}$-order spherical harmonics by projections $\Pi^{(n)}_X$  
and $\Pi^{(\mathbf s)}_K$, where compound index $\mathbf s=(s,f,a)$
indicates the degrees both of the target and original functions. Like
the projections introduced regarding the X-variable, also the latter
projections are the polynomials of the Laplacian $\Delta_{K_u}$.
More precisely, they appear in the form 
$q_{\mathbf s}\Pi^{(s)}_K\Delta^{(f+a-s)/2}_{K_u}$,
where the first term is a constant and the last term produces from
a $(f+a)^{th}$-order polynomial of the K-variable an 
$s^{th}$-order polynomial which is then projected to the $s^{th}$-order
spherical harmonics in the same way how it was defined regarding the
X-variable.
(These computations are described in Section 6.4 of \cite{sz8}.)

The corresponding formula for constructing functions by which the
boundary conditions can be handled is:
\begin{equation}
\label{newDiracwave2}
\int_{\mathbb R^l}\sum_{\mathbf s}\phi_{\mathbf s}(|X|,|K|)
\Pi^{\mathbf s}_K(\varphi (K_u)\Pi^{(n)}_X(\prod z_i^{p_i}
\overline z_i^{q_i})(K_u,X))
e^{\mathbf i\langle K,Z\rangle}dK.
\end{equation}
It is also pointed out in \cite{sz8} that for fixed spherical harmonics 
$\varphi (K_u)$, non-trivial projections are defined just for those 
$s$-values which satisfy the inequality 
$(f-a)\leq s\leq (f+a)$ and $s$ has the same parity as  
$(f-a)$ resp. $(f+a)$. Furthermore,
functions defined by different $a$'s are
projected into independent subspaces.   
That is, this formula generically involves $(a+1)$-tuples,
$(\phi_{(f-a)},\dots ,\phi_{(f+a)})$, of functions 
depending $|X|$ and $|K|$.

The functions constructed in this way form a larger space, 
$\mathbb F_Z(\mathbb {LT}w_{\mathbf B}^{(p_i,q_i)})$,
then those
constructed in (\ref{newDiracwave}). Since functions (\ref{newDiracwave}) 
form an everywhere dense subspace in the whole straightly defined space, 
they are everywhere dense also in space spanned by functions defined in
(\ref{newDiracwave2}). This means that in the newly defined space the
functions can just be approximated by the previous functions. 
However, no natural (obvious) approximation exist what would make
the consideration of the new formula evident. 
Anyhow, the statement about the density of these function spaces is
enough to establish the following statement: {\it For any two members 
of an isospectrality family, the intertwining operator
defined for the primarily functions extends to the newly defined ones 
such that it 
associates functions defined in terms of same $(X,Z)$-radial functions
and the corresponding projected functions.} More precisely, this is the
only option for a continuous extension, Indeed,
if the operator is defined in this way, then on
$\mathbb F_Z(\mathbb {T}w_{\mathbf B}^{(p_i,q_i)})$ it is the same  
as the operator defined originally. This argument implies
that the operator defined in terms of the newly defined functions
is still an operator intertwining the Laplacians.

For investigating the boundary conditions there are two important tools
established. One of them is the Hankel transform, which is proved in
Section 6.3 of \cite{sz8} in the following form:  
{\it Any $s^{th}$-order spherical harmonic
$\zeta_l^{(s)}(K_u)$ defined on the unit sphere of $\mathbb R^l$
defines, by the formula 
\begin{equation}
\mathcal H_l^{(s)}(\phi)(|Z|)\zeta_l^{(s)}(Z_u)=\int_{\mathbb R^l}
\phi (|K|)\zeta_l^{(s)}(K_u)e^{\mathbf i\langle Z,K\rangle}dK,
\end{equation}
a uniquely determined transformation
$\mathcal H_l^{(s)}(\phi)(|Z|)$  on the $L^2$-radial
functions $\phi (|Z|)$ which depends just on $s$ and $l$.}
This statement implies that the newly defined
function space is appropriate for constructing the complete space 
both of the Dirichlet and Z-Neumann functions 
in terms of the $(X,Z)$-radial
functions. It turns out that these two conditions can be characterized
as being scalar, meaning that they are satisfied if and only if functions 
$\mathcal H_l^{(s)}(\phi_{\mathbf s})(|X|,|Z|)$ appearing in the above 
formulas satisfy them individually. Actually, these functions are also
explicitly described in \cite{sz8}, establishing this 
statement a much stronger form. As a result, these functions 
are intertwined by the above operator, indeed. 

The other important
tool is the {\it inner algorithm} (cf. Section 6.7 in \cite{sz8}) 
induced by the action of the
angular momentum operator or operator $D_Z\bullet$ on the second
Z-Fourier transform formula. In this algorithm this action is iterated
such that in the $k^{th}$-step the desired functions are constructed by
functions $\phi^{(k-1)}_{\mathbf s}$ obtained in the previous step
by using Hankel transform, radial 
differentiation, and certain combinations of them which can be 
described as averaging by the roulette operator \cite{sz8}. The point is 
that this process involves just the radial functions and it ends up 
either in finite or infinite steps. In the latter case the 
sought functions are obtained by the limiting $k\to\infty$.  

There is far more difficult problem
to establish the above statement also for functions yielding 
the regular Neumann condition. The source of these difficulties is
that, contrary to the Dirichlet and Z-Neumann conditions, 
this one does not break down to single individual functions, 
but, it can be expressed in terms of all 
functions $(\phi_{(f-a)},\dots ,\phi_{(f+a)})$.
More precisely, with the help of the
inner algorithm, one can construct a compound Neumann operator
$\mathbb N^{(f,a)}_{s}(\phi_{(f-a)},\dots ,\phi_{(f+a)})$ acting
on radial functions such that a function constructed by
(\ref{newDiracwave2}) satisfies the regular Neumann condition if and
only if the radial functions defined by the compound Neumann operator 
vanishes at the boundary points. That is, also this condition can be 
expressed in terms of radial
functions, therefore it is also intertwined by the above operator.

The complexity of angular momentum operator $\mathbf M_Z$ is fascinating.
Its action can be described by the Hankel transform and the above mentioned
inner algorithm in a more precise way (cf. Sections 6.7 and 6.8
of \cite{sz8}). These tools reveal that it appears as the 
sum of an extrinsic orbiting operator, $\mathbf{L}$, operating on the
radial functions $\phi_{\mathbf s}$ without defining
permutations (i. e., averaging by the roulette operator) on them, 
and an intrinsic spin operator, $\mathbf{S}$, 
acting on radial functions, but contrary to 
$\mathbf{L}$, it defines also an irreducible permutation on 
functions $\phi_s$. The actual appearance of this operator on the above
functions is:
\begin{equation}
\label{newDiracwave3}
\int_{\mathbb R^l}\bigcirc_{\mathbf s}(\phi_{(f-a)},\dots ,\phi_{(f+a)})
\Pi^{(\mathbf s)}_K(\varphi \Pi^{(n)}_X(\prod z_i^{p_i}
\overline z_i^{q_i})(K_u,X))
e^{\mathbf i\langle K,Z\rangle}dK.
\end{equation}
The intrinsic life of particles is encoded into $\mathbf{S}$. 
Also the strong nuclear forces, keeping the 
particles having interior together, can be explained by this operator.
 
The rest part,
$\OE =\Delta_X+(1+\frac 1{4}|X|^2)\Delta_Z+\mathbf{L}$,
of the complete Laplacian $\Delta$ exhibits just orbiting spin. Its
action can be reduced to a single radial function. More precisely,  
the exterior operator $\OE$, 
on constant radius Z-ball bundles 
reduces to a radial operator appearing in terms of the 
Dirichlet-, Neumann-, resp. mixed-condition-eigenvalues 
$\lambda_i^{(s)}$ of the 
Z-ball $B_Z(R)$ in the form 
\begin{equation}
\label{BLf_mu}
(\Diamond_{\mu ,\tilde t}\mathtt f)(\tilde t)=
4\tilde t\mathtt f^{\prime\prime}(\tilde t)
+(2k+4n)\mathtt f^\prime (\tilde t)
-(2m{\mu} +4{\mu^2} (1 +
{1\over 4}\tilde t))\mathtt f(\tilde t).
\end{equation}
By the substitution 
$\mu =\sqrt{\lambda_i^{(s)}/4}$, this is exactly
the radial Ginsburg-Landau-Zeeman operator 
(\ref{Lf_lambda}) obtained on Z-crystal models. 

The physical forces corresponding to operator $\OE$ are the weak nuclear
forces by which the beta decays are explained. The theoretical 
establishment of this force went through enormous developments 
\cite{v,w1,w2}.
It started out with Fermi's theory which was highly surpassed by
the Glashow-Weinberg-Salam theory, whose greatest achievement was the 
unification of the weak nuclear force with the electromagnetic force.
The above statement, which is unknown both in physics and mathematics, 
is an exact mathematical establishment of this unification.
In physics the unification with the other forces, that is with strong 
nuclear forces and gravitation, is one of the most 
intensely investigated unsolved questions. Since 
the Z-crystal and Z-ball-bundle models unify the 
electromagnetic and weak nuclear forces also with the strong nuclear
forces, this unification is actually much stronger than those 
established in the Glashow-Weinberg-Salam theory. 
The main unifying idea is that all these forces can be derived 
from the very same operator $\Delta$.
They are distinguished only by the invariant subspaces to which
the $\Delta$ is restricted to. That is, the individual forces
emerge on these individual invariant subspaces separately.   

\noindent{\bf Conclusion:} {\it On these mathematical-physical 
models the investigated isospectralities
are equivalent to the spectral C/I-symmetry known in
elementary particle physics. If the elements of the basis $\mathbf B$
are picked up from the invariant subspaces $\mathcal X^{(a)}$ and
$\mathcal X^{(b)}$, then the isospectrality is the manifestation of
the pure spectral C-symmetry. 
In other words, the isospectrality proofs in my constructions 
mathematically demonstrate that, regarding the spectrum of
the Hamilton operator, these physical models obey the 
physical C-symmetry law. That is, this spectrum is not changing if some
of the particles are exchanged for their antiparticles.}

\subsection{Wave mechanics.} 
The wave operators corresponding to the above
Hamilton operators emerge in the Laplacian of the static resp. 
solvable extensions of the nilpotent groups. In physics, where this
operator should be a hyperbolic wave operator, the 
extended metric must be indefinite having Lorenz signature, where the
time axis has signature $-1$. Regarding both bundles, the solutions of
the wave equations considered on the whole bundle can be written up 
by explicit integral formulas (cf. Section 7 of \cite{sz8}). 
On ball-bundles,
corresponding to the static and solvable extensions, 
these de Broglie wave packets appear in the form
\begin{eqnarray}
\label{newDiracwave3}
\int_{\mathbb R^l}\phi_s
\Pi^{(s)}_K(\varphi \Pi^{(n)}_X(\prod z_i^{p_i}
\overline z_i^{q_i}))
e^{\mathbf i(\langle K,Z\rangle -\omega t)}dK,\quad {\rm resp.}
\\
\int_{\mathbb R^l}\phi_s
\Pi^{(s)}_K(\varphi \Pi^{(n)}_X(\prod z_i^{p_i}
\overline z_i^{q_i}))
e^{\mathbf i(\langle K,Z\rangle -\omega e^T)}dK,
\nonumber
\end{eqnarray}
where $\omega$ is a constant depending on the particle system. On 
Z-crystals, an appropriate discrete version of the integral defines
these wave packets. Then, on the extended Z-crystals, 
the Laplacian is the sum
of a Schr\"odinger operator determined for electron-positron systems
and an electron-positron-neutrino operator accompanying the
electron-positron system. On ball-bundles, the corresponding
particles regarding the $\OE$-operator are the W- and Z-particles, 
introduced by Weinberg in his weak-force theory. Analogous particles can be
introduced regarding the complete Laplacian $\OE +\mathbf S$. 
These details are completely described in Section 7 of \cite{sz8}.

On manifolds endowed with indefinite metrics the isospectrality questions
can not be raised in general. However, this question can be raised
regarding positive definite extensions, where one extends 
the natural positive definite
inner product defined at $(0,0,1)$ on the tangent space. These manifolds
can be considered as non-relativistic Newtonian space-time models.  
Actually the harmonic
isospectral manifolds discussed in the AriasMarco-Sch\"uth paper
appear exactly among these manifolds. The metrics on nilpotent groups
are neither harmonic nor Einstein thus they could not have 
been effected by the AM\&S-theorem even if it was a right statement. 
The proof of isospectrality on the positive definite
solvable extensions can be established similarly as on the nilpotent
groups. In this case the radial functions should be considered in the
form $\phi_s(|X|,|K|,t)$, where $t>0$ is the new parameter added to 
$(X,Z)$ which indicates that the solvable extensions are defined by 
the half-space extension of two-step nilpotent Lie groups. These
details and the isospectralities on the boundaries of Z-ball-bundles are 
described in \cite{sz6}.
   
\medskip\medskip

\noindent{\bf The problems arising in my papers \cite{sz4,sz5}.}
The intertwining operators are not correctly established in my papers
cited above. These difficulties concern only these papers, where the
isospectrality examples are constructed on Z-ball resp. Z-sphere bundles.
They do not effect my other constructions performed on Z-torus bundles, 
and, even in these papers, they concern only the construction of
the intertwining operators, while the other statements are unaffected. 

These operators are constructed in \cite{sz4,sz5} by means
of a fixed complex 
structure $J_0\in E_{skew}(\mathcal X)$. 
The complex structure obtained by $\sigma$ deformation is denoted
by $J_0^\prime$. Then, the intertwining operator, $\kappa$, is defined 
in these papers by the correspondence 

\begin{eqnarray}
\label{oldkappa}
\kappa :\quad A(|X|,Z)\Pi_X^{(n)}(\prod_i \langle Q_i+
\mathbf iJ_{0}(Q_i),X\rangle^{p_i}
\langle Q_i-\mathbf iJ_{0}(Q_i),X\rangle^{q_i})
\\
\to\quad
A(|X|,Z)\Pi_X^{(n)}(\prod_i \langle Q_i+\mathbf iJ^\prime_{0}(Q_i),X
\rangle^{p_i}
\langle Q_i-\mathbf iJ^\prime_{0}(Q_i),X\rangle^{q_i})
\nonumber
\end{eqnarray}
where $Q_i\in\mathcal X$ are arbitrary X-vectors. That is, the associated
functions appear in terms of the associated complex structures 
$J_0$ resp. $J_0^\prime$ in the very same form. It obviously intertwines
these two complex structures 

It was H. F\"urstenau \cite{f} who recognized that this map was not 
well defined. For if one considers a fixed complex
basis $\mathbf B=\{B_1,\dots,B_{k/2}\}$ on the X-space, then the above 
correspondence can be prescribed only for functions defined by 
such $Q_i$'s which are in 
the real subspace $Span_{\mathbb R}(\mathbf B)$ spanned by real linear 
combinations of the basis-elements $B_i$. It is evident, that these 
correspondences will define the
action of $\kappa$ on the rest of functions defined by other 
$Q_i$'s not lying
in $Span_{\mathbb R}(\mathbf B)$. The very same map can be defined
in terms of the complex valued polynomials constructed above regarding the
basis $\mathbf B$. That is, it is enough to consider only the 
correspondence which associate polynomials (appearing in the same form
regarding the associated complex structures) to each other. 
It is also clear that the action of this well defined $\kappa$ on
functions defined by $Q_i$'s not lying in 
$Span_{\mathbb R}(\mathbf B)$ differs from the action described in 
(\ref{oldkappa}). Thus the above map
is not well defined, indeed. (F\"urstenau used different arguments for
explaining this problem.) 

This recognition triggered the reconstruction of the above ill-defined
intertwining operator.
Further investigations showed that a correct operator could have not been
established by a single complex structure $J_0$. Instead, all
complex structures $J_{Z_u}$ must be involved to its definition. 
Even the exchange
of $J_0$ for $J_{Z_u}$ does not alter this newly defined operator 
into a correct one. It becomes the desired intertwining operator,
however, if one applies also the Z-Fourier transform to the latter 
altered functions.
This step is really necessary because the operator defined in terms
of functions standing behind the integral sign intertwines just those
physical operators which appear after the Laplacian enters behind the
integral sign.
In terms of a complex basis $\mathbf B$ and polynomials written up 
regarding this basis, the reconstructed intertwining operator can 
be defined by the functions introduced by the twisted Z-Fourier transforms 
(\ref{newDiracwave}) and (\ref{newDiracwave2}).

These reconstructed formulas were described first in \cite{sz6} and 
lectured about them also at the international conference held at CUNY, in 
February of 2006 \cite{sz7}. The physical content reviewed in this paper
is the result of a recent development. It is described in \cite{sz8}. 
In this paper a completely new mathematical idea, namely the Hankel 
transform, is introduced, by which both the electroweak and 
strong interaction Hamilton operators can be established. By this
transform simple transparent proofs can be given also for 
mathematical theorems such as the density theorems and several other
statements concerning intertwining of Laplacians or boundary 
conditions. In \cite{sz6} these statements
are established by a different integral transform, called dual
Radon transform. An other statement used there is the so 
called independence theorem. It should be pointed out that
the exploration of the physical contents inherent in these
structures is not quite much efficient by these tools than by the Hankel 
transform. This statement is certainly true for electroweak and strong 
interactions, which become understandable just by the Hankel transform. 
The strength of the latter method can be demonstrated also by the fact 
that all theorems established in \cite{sz6} can be established also by 
this tool.

\end{document}